\documentclass[reqno,11pt]{amsart}
\usepackage{amsfonts}

\usepackage{graphicx}
\usepackage{pstricks}
\usepackage{amsmath}
\usepackage{amsmath}
\usepackage{amsxtra}

\theoremstyle{plain}
\newtheorem{theorem}{Theorem}[section]
\newtheorem{lemma}[theorem]{Lemma}

\newtheorem{proposition}[theorem]{Proposition}
\newtheorem{corollary}[theorem]{Corollary}

\newtheorem{definition}[theorem]{Definition}

\theoremstyle{definition}

\newtheorem{example}[theorem]{Example}
\newtheorem{remark}[theorem]{Remark}

\numberwithin{equation}{section}
\begin{document}
\title{Propagation phenomena\\
for hyponormal 2-variable weighted shifts}
\author{Ra\'{u}l E. Curto}
\address{Department of Mathematics, The University of Iowa, Iowa City, Iowa
52242}
\email{rcurto@math.uiowa.edu}
\urladdr{http://www.math.uiowa.edu/\symbol{126}rcurto/}
\author{Jasang Yoon}
\address{Department of Mathematics, Iowa State University, Ames, Iowa 50011}
\email{jyoon@iastate.edu}
\urladdr{http://www.public.iastate.edu/\symbol{126}jyoon/}
\thanks{The first named author was partially supported by NSF Grants
DMS-0099357 and DMS-0400471.}

\begin{abstract}
We study the class of hyponormal $2$-variable weighted shifts with two
consecutive equal weights in the weight sequence of one of the coordinate
operators. We show that under natural assumptions on the coordinate
operators, the presence of consecutive equal weights leads to horizontal or
vertical flatness, in a way that resembles the situation for $1$-variable
weighted shifts. \ In $1$-variable, it is well known that flat weighted
shifts are necessarily subnormal (with finitely atomic Berger measures). \
By contrast, we exhibit a large collection of flat (i.e., horizontally and
vertically flat) $2$-variable weighted shifts which are hyponormal but not
subnormal. \ Moreover we completely characterize the hyponormality and
subnormality of symmetrically flat contractive $2$-variable weighted shifts.
\end{abstract}

\maketitle

\section{\label{Int}Introduction}

The Lifting Problem for Commuting Subnormals (LPCS) asks for necessary and
sufficient conditions on a pair of commuting subnormal operators on Hilbert
space to admit a joint normal extension. \ In previous work we have proved
that the (joint) hyponormality of the pair, while being a necessary
condition, is by no means sufficient (\cite{CuYo1}, \cite{CuYo2}). \ We have
also established that in a very special situation, hyponormality is indeed
sufficient \cite[Theorem 5.2 and Remark 5.3]{CuYo1}. \ This involves $2$%
-variable weighted shifts with weight sequences which are constant except
for the $0$-th row in the index set $\mathbb{Z}_{+}^{2}$. \ One is then
tempted to claim that a similar result might be true for weight sequences
which are constant in a slightly smaller domain of indices, e.g., those
indices $\mathbf{k}\in \mathbb{Z}_{+}^{2}$ with $k_{1},k_{2}\geq 1$. \
However, in this paper we show that such is not the case, that is,
hyponormality and subnormality are quite different even in those cases.

For $\alpha \equiv \{\alpha _{k}\}_{k=0}^{\infty }$ a bounded sequence of
positive real numbers (called \textit{weights}), let $W_{\alpha }:\ell ^{2}(%
\mathbb{Z}_{+})\rightarrow \ell ^{2}(\mathbb{Z}_{+})$ be the associated 
\textit{unilateral weighted shift}, defined by $W_{\alpha }e_{k}:=\alpha
_{k}e_{k+1}\;($all $k\geq 0)$, where $\{e_{k}\}_{k=0}^{\infty }$ is the
canonical orthonormal basis in $\ell ^{2}(\mathbb{Z}_{+}).$ \ A
quadratically hyponormal weighted shift $W_{\alpha }$ with $\alpha
_{k+1}=\alpha _{k}$ for some $k\geq 1$ must necessarily be (i) \textit{flat}
(i.e., $\alpha _{1}=\alpha _{2}=\alpha _{3}=\cdots $), and (ii) subnormal. \
For $2$-variable weighted shifts associated with weight sequences $\{\alpha
_{\mathbf{k}}\},\{\beta _{\mathbf{k}}\}\in \ell ^{\infty }(\mathbb{Z}%
_{+}^{2})$, we first establish the correct analogue of (i) (Theorem \ref%
{thm33}), and we then show that there is a rich family of sequences $%
\{\alpha _{\mathbf{k}}\},\{\beta _{\mathbf{k}}\}$ giving rise to flat,
non-subnormal, hyponormal $2$-variable weighted shifts; this is in sharp
contrast with the $1$-variable situation. \ The optimality of Theorem \ref%
{thm33} is established through an elaborate construction which uses
Bergman-like weighted shifts (Theorem \ref{Propagation of Induction}). \
Finally, in Section \ref{symm} we completely characterize the hyponormality
and subnormality of symmetrically flat contractive $2$-variable weighted
shifts, which sheds new light on the relationship between flatness and
subnormality.

Recall that a bounded linear operator $T\in \mathcal{B}(\mathcal{H})$ on a
complex Hilbert space $\mathcal{H}$ is \textit{normal} if $T^{\ast
}T=TT^{\ast }$, \textit{subnormal} if $T=N|_{\mathcal{H}}$, where $N$ is
normal and $N(\mathcal{H})\subseteq \mathcal{H}$, and \textit{hyponormal} if 
$T^{\ast }T\geq TT^{\ast }$.\ \ For $k\geq 1$ and $T\in \mathcal{B}(\mathcal{%
H}),$ $T$ is $k$\textit{-hyponormal} if $(I,T,\cdots ,T^{k})$ is (jointly)
hyponormal. \ Additionally, $T$ is \textit{weakly }$\mathit{k}$\textit{%
-hyponormal} if $p(T)$ is hyponormal for every polynomial $p$ of degree at
most $k$. \ Thus, if $T$ is $k$-hyponormal then $T$ is weakly $k$%
-hyponormal, and ``hyponormality'' ``$1$-hyponormality''\ and ``weak $1$%
-hyponormality''\ are all identical notions (\cite{Ath}). \ On the other
hand, results in (\cite{CMX}), (\cite{QHWS}) and (\cite{McCP}) show that if $%
T$ is weakly $2$-hyponormal (also called \textit{quadratically hyponormal}),
then $T$ need not be $2$-hyponormal. \ The Bram-Halmos characterization of
subnormality (\cite[III.1.9]{Con}) can be paraphrased as follow: \ $T$ is
subnormal if and only if $T$ is $k$-hyponormal for every $k\geq 1$ (%
\cite[Proposition 1.9]{CMX}). \ In particular, each subnormal operator is 
\textit{polynomially hyponormal} (i.e., weakly $k$-hyponormal for every $%
k\geq 1$). \ The converse implication, whether $T$ polynomially hyponormal $%
\Rightarrow $ $T$ subnormal, was settled in the negative in (\cite{CuPu});
indeed, it was shown that there exists a polynomially hyponormal operator
which is not $2$-hyponormal. \ Previously, S. McCullough and V. Paulsen had
established (\cite{McCP}) that one can find a non-subnormal polynomially
hyponormal operator if and only if one can find a unilateral weighted shift
with the same property. \ Thus, although the existence proof in (\cite{CuPu}%
) is abstract, by combining the results in (\cite{CuPu}) and (\cite{McCP})
we now know that there exists a polynomially hyponormal unilateral weighted
shift which is not subnormal.

For $S,T\in \mathcal{B}(\mathcal{H})$ we let $[S,T]:=ST-TS$. \ We say that a
commuting $n$-tuple $\mathbf{T}=(T_{1},\cdots ,T_{n})$ of operators on $%
\mathcal{H}$ is (jointly) \textit{hyponormal} if the operator matrix 
\begin{equation*}
\lbrack \mathbf{T}^{\ast },\mathbf{T]:=}\left( 
\begin{array}{llll}
\lbrack T_{1}^{\ast },T_{1}] & [T_{2}^{\ast },T_{1}] & \cdots & [T_{n}^{\ast
},T_{1}] \\ 
\lbrack T_{1}^{\ast },T_{2}] & [T_{2}^{\ast },T_{2}] & \cdots & [T_{n}^{\ast
},T_{2}] \\ 
\text{ \thinspace \thinspace \quad }\vdots & \text{ \thinspace \thinspace
\quad }\vdots & \ddots & \text{ \thinspace \thinspace \quad }\vdots \\ 
\lbrack T_{1}^{\ast },T_{n}] & [T_{2}^{\ast },T_{n}] & \cdots & [T_{n}^{\ast
},T_{n}]%
\end{array}%
\right)
\end{equation*}%
is positive on the direct sum of $n$ copies of $\mathcal{H}$ (cf. \cite{Ath}%
, \cite{CMX}). \ The $n$-tuple $\mathbf{T}$ is said to be \textit{normal} if 
$\mathbf{T}$ is commuting and each $T_{i}$ is normal, and $\mathbf{T}$ is 
\textit{subnormal }if $\mathbf{T}$ is the restriction of a normal $n$-tuple
to a common invariant subspace. \ Clearly, normal $\Rightarrow $ subnormal $%
\Rightarrow $ hyponormal.

For $\alpha \equiv \{\alpha _{k}\}_{k=0}^{\infty }\in \ell ^{\infty }(%
\mathbb{Z}_{+})$ and $W_{\alpha }$ the associated unilateral weighted shift,
the \textit{moments} of $\alpha $ are given as 
\begin{equation*}
\gamma _{k}\equiv \gamma _{k}(\alpha ):=\left\{ 
\begin{array}{cc}
1 & \text{if }k=0 \\ 
\alpha _{0}^{2}\cdot ...\cdot \alpha _{k-1}^{2} & \text{if }k>0%
\end{array}%
\right. .
\end{equation*}%
It is easy to see that $W_{\alpha }$ is never normal, and that it is
hyponormal if and only if $\alpha _{0}\leq \alpha _{1}\leq \cdots $. \ If $%
\alpha _{k+1}=\alpha _{k}$ for all $k\geq 1,$ $W_{\alpha }$ is called 
\textit{flat}. \ On occasion, we will write $shift(\alpha _{0},\alpha
_{1},\alpha _{2},\cdots )$ to denote the weighted shift with weight sequence 
$\{\alpha _{k}\}_{k=0}^{\infty }$. \ We also denote by $U_{+}:=shift(1,1,1,%
\cdots )$ the (unweighted) unilateral shift, and for $0<a<1$ we let $%
S_{a}:=shift(a,1,1,\cdots )$; the shift $S_{a}$ is the prototypical flat
weighted shift, and it is subnormal.

Similarly, consider double-indexed positive bounded sequences $\{\alpha _{%
\mathbf{k}}\},\{\beta _{\mathbf{k}}\}\in \ell ^{\infty }(\mathbb{Z}_{+}^{2})$
, $\mathbf{k}\equiv (k_{1},k_{2})\in \mathbb{Z}_{+}^{2}:=\mathbb{Z}%
_{+}\times \mathbb{Z}_{+}$ and let $\ell ^{2}(\mathbb{Z}_{+}^{2})$\ be the
Hilbert space of square-summable complex sequences indexed by $\mathbb{Z}%
_{+}^{2}$. \ (Recall that $\ell ^{2}(\mathbb{Z}_{+}^{2})$ is canonically
isometrically isomorphic to $\ell ^{2}(\mathbb{Z}_{+})\bigotimes \ell ^{2}(%
\mathbb{Z}_{+})$.) \ We define the $2$\textit{-variable weighted shift} $%
\mathbf{T}$\ by 
\begin{equation*}
T_{1}e_{\mathbf{k}}:=\alpha _{\mathbf{k}}e_{\mathbf{k+\varepsilon }_{1}}
\end{equation*}%
\begin{equation*}
T_{2}e_{\mathbf{k}}:=\beta _{\mathbf{k}}e_{\mathbf{k+\varepsilon }_{2}},
\end{equation*}%
where \textbf{$\varepsilon $}$_{1}:=(1,0)$ and \textbf{$\varepsilon $}$%
_{2}:=(0,1)$. \ Clearly, 
\begin{equation}
T_{1}T_{2}=T_{2}T_{1}\Longleftrightarrow \beta _{\mathbf{k+\varepsilon }%
_{1}}\alpha _{\mathbf{k}}=\alpha _{\mathbf{k+\varepsilon }_{2}}\beta _{%
\mathbf{k}}\;\;(\text{all }\mathbf{k}).  \label{commuting}
\end{equation}%
In an entirely similar way one can define multivariable weighted shifts. \
Trivially, a pair of unilateral weighted shifts $W_{\alpha }$ and $W_{\beta
} $ gives rise to a $2$-variable weighted shift $\mathbf{T}\equiv
(T_{1},T_{2}) $, if we let $\alpha _{(k_{1},k_{2})}:=\alpha _{k_{1}}$ and $%
\beta _{(k_{1},k_{2})}:=\beta _{k_{2}}\;($all $k_{1},k_{2}\in \mathbb{Z}%
_{+}^{2})$. \ In this case, $\mathbf{T}$ is subnormal (resp. hyponormal) if
and only if so are $T_{1}$ and $T_{2}$; in fact, under the canonical
identification of $\ell ^{2}(\mathbb{Z}_{+}^{2})$ and $\ell ^{2}(\mathbb{Z}%
_{+})\bigotimes \ell ^{2}(\mathbb{Z}_{+})$, $T_{1}\cong W_{\alpha
}\bigotimes I$ and $T_{2}\cong I\bigotimes W_{\beta }$, and $\mathbf{T}$ is
also doubly commuting. \ For this reason, we do not focus attention on
shifts of this type, and use them only when the above mentioned triviality
is desirable or needed.

We now recall a well known characterization of subnormality for single
variable weighted shifts, due to C. Berger (cf. \cite[III.8.16]{Con}): \ $%
W_{\alpha }$ is subnormal if and only if there exists a probability measure $%
\xi $ supported in $[0,\left\| W_{\alpha }\right\| ^{2}]$ (called the 
\textit{Berger measure} of $W_{\alpha }$) such that $\gamma _{k}(\alpha
):=\alpha _{0}^{2}\cdot ...\cdot \alpha _{k-1}^{2}=\int t^{k}d\xi (t)$ $%
(k\geq 1)$. \ For instance, the Berger measures of $U_{+}$ and $S_{a}$ are $%
\delta _{1}$ and $(1-a^{2})\delta _{0}+a^{2}\delta _{1}$, respectively,
where $\delta _{x}$ denotes the point-mass probability measure with support
the singleton $\{x\}$.

If $W_{\alpha }$ is subnormal, and if for $h\geq 1$ we let $\mathcal{M}%
_{h}:=\bigvee \{e_{k}:k\geq h\}$ denote the invariant subspace obtained by
removing the first $h$ vectors in the canonical orthonormal basis of $\ell
^{2}(\mathbb{Z}_{+})$, then the Berger measure of $W_{\alpha }|_{\mathcal{M}%
_{h}}$ is $\frac{1}{\gamma _{h}}t^{h}d\xi (t)$. \ For $h=2$, one can use
this to prove the following result.

\begin{lemma}
\label{lembeta1}Let $T\equiv shift(\beta _{0},\beta _{1},\cdots )$ be a
subnormal weighted shift, with Berger measure $\eta $, and let $T_{\mathcal{M%
}}$ be its restriction to $\mathcal{M}:=\vee \{e_{2},e_{3},\cdots \}$. \
Then $\beta _{1}^{2}=(\left\| \frac{1}{t}\right\| _{L^{1}(\eta _{\mathcal{M}%
})})^{-1}$.
\end{lemma}

\begin{proof}
We have 
\begin{equation*}
\left\| \frac{1}{t}\right\| _{L^{1}(\eta _{\mathcal{M}})}=\int \frac{1}{t}(%
\frac{1}{\gamma _{2}}t^{2}d\eta (t))=\frac{1}{\gamma _{2}}\int td\eta (t)=%
\frac{\gamma _{1}}{\gamma _{2}}=\frac{1}{\beta _{1}^{2}},
\end{equation*}%
as desired.
\end{proof}

\begin{corollary}
\label{corbeta1}Let $\mathbf{T}\equiv (T_{1},T_{2})$ be a commuting $2$%
-variable weighted shift, assume that $T_{2}$ is subnormal, and assume that
there exists $k_{1}\geq 0$ such that $\alpha _{(k_{1},k_{2})}=\alpha
_{(k_{1},k_{2})+\mathbf{\varepsilon }_{2}}$ for all $k_{2}\geq 2$. \ Then $%
\beta _{(k_{1},1)}=\beta _{(k_{1},1)+\mathbf{\varepsilon }_{1}}$.
\end{corollary}

\begin{proof}
Consider $S:=shift(\beta _{(k_{1},2)},\beta _{(k_{1},3)},\cdots )$ and $%
S^{\prime }:=shift(\beta _{(k_{1+1},2)},\beta _{(k_{1}+1,3)},\cdots )$. \
Since $T_{2}$ is subnormal, we know that both $S$ and $S^{\prime }$ are
subnormal, with Berger measures $\eta $ and $\eta ^{\prime }$, respectively.
\ Since $\alpha _{(k_{1},k_{2})}=\alpha _{(k_{1},k_{2})+\mathbf{\varepsilon }%
_{2}}$, the commuting property (\ref{commuting}) readily implies that $\beta
_{(k_{1},k_{2})}=\beta _{(k_{1},k_{2})+\mathbf{\varepsilon }_{1}}$ for all $%
k_{2}\geq 2$, that is $S=S^{\prime }$, that is, $\eta =\eta ^{\prime }$. \
By Lemma \ref{lembeta1}, we must have 
\begin{equation*}
\beta _{(k_{1},1)}^{2}=(\left\Vert \frac{1}{t}\right\Vert _{L^{1}(\eta
)})^{-1}=(\left\Vert \frac{1}{t}\right\Vert _{L^{1}(\eta ^{\prime
})})^{-1}=\beta _{(k_{1},1)+\mathbf{\varepsilon }_{1}}^{2},
\end{equation*}%
as desired.
\end{proof}

\textit{Acknowledgments}. \ The authors are grateful to the referee for
several suggestions which helped improved the presentation. \ Example \ref%
{numericalex1} and some of the proofs of the results in this paper were
obtained using calculations with the software tool \textit{Mathematica \cite%
{Wol}.}

\section{\label{results}Propagation Phenomena for $1$-Variable Weighted
Shifts}

In this section, we review some basic propagation phenomena for $1$-variable
weighted shifts, and we then develop the results for the $2$-variable case
in Sections \ref{Section3} and \ref{Section4}. \ J. Stampfli showed in \cite%
{Sta} that for a subnormal weighted shift $W_{\alpha }$, a propagation
phenomenon occurs which forces the flatness of $W_{\alpha }$ whenever two
equal weights are present.

\begin{proposition}
\label{stasub}(Subnormality, One-variable Case) (\cite{Sta}) \ Let $%
W_{\alpha }$ be a subnormal weighted shift with weight sequence $\{\alpha
_{k}\}_{k=0}^{\infty }.$ \ If $\alpha _{k}=\alpha _{k+1}$ for some $k\geq 0$%
, then $W_{\alpha }$ is flat.
\end{proposition}

The first named author showed that in the presence of $2$-hyponormality
(resp. quadratic hyponormality) of weighted shifts, a propagation phenomenon
also occurs which forces the flatness of $W_{\alpha }$ whenever two equal
weights (resp. three equal weights) are present.

\begin{proposition}
($2$-hyponormality, One-variable Case) (\cite{QHWS}) \ Let $W_{\alpha }$ be
a $2$-hyponormal weighted shift with weight sequence $\{\alpha
_{k}\}_{k=0}^{\infty }.$ \ If $\alpha _{k}=\alpha _{k+1}$ for some $k\geq 0$%
, then $W_{\alpha }$ is flat.
\end{proposition}

\begin{proposition}
\label{propCu}(Quadratic Hyponormality, One-variable Case) (\cite{QHWS}) \
Let $W_{\alpha }$ be a unilateral weighted shift with weight sequence $%
\{\alpha _{k}\}_{k=0}^{\infty },$ and assume that $W_{\alpha }$ is
quadratically hyponormal. $\ $If $\alpha _{k}=\alpha _{k+1}=$ $\alpha _{k+2}$
for some $k\geq 0$, then $W_{\alpha }$ is flat.
\end{proposition}

Y. Choi later improved Proposition \ref{propCu}, as follows.

\begin{proposition}
\label{quadhypo3}(Quadratic Hyponormality, One-variable Case, Improved
Version) (\cite{Ch}) \ Let $W_{\alpha }$ be a unilateral weighted shift with
weight sequence $\{\alpha _{k}\}_{k=0}^{\infty },$ and assume that $%
W_{\alpha }$ is quadratically hyponormal. $\ $If $\alpha _{k}=\alpha _{k+1}$
for some $k\geq 1$, then $W_{\alpha }$ is flat.
\end{proposition}

Moreover, Y. Choi showed that, in the presence of polynomially
hyponormality, two consecutive equal weights again force flatness.

\begin{proposition}
(Polynomially hyponormality) (\cite{Ch}) \ Let $W_{\alpha }$ be a unilateral
weighted shift with weight sequence $\{\alpha _{k}\}_{k=0}^{\infty },$ and
assume that $W_{\alpha }$ is polynomially hyponormal. $\ $If $\alpha
_{k}=\alpha _{k+1}$ for some $k\geq 0$, then $W_{\alpha }$ is flat.
\end{proposition}

\section{\label{Section3}Propagation in the $2$-variable Hyponormal Case}

In this section, we show that if a commuting, (jointly) hyponormal pair $%
\mathbf{T\equiv }(T_{1},T_{2})$ with $T_{1}$ quadratically hyponormal
satisfies $\alpha _{(k_{1}+1,k_{2})}=\alpha _{(k_{1},k_{2})}$ for some $%
k_{1},k_{2}\geq 1$, then $(T_{1},T_{2}(I\otimes U_{+}^{k_{2}-1}))$ is
horizontally flat (see Definition\ \ref{def31} below); this is the content
of Theorem \ref{thm33}. \ We also prove that Theorem \ref{thm33} is optimal
in the following sense: the propagation does not extend either to the left ($%
0$-th column) or down (below $k_{2}$-th level).

We begin with:

\begin{definition}
\label{def31}A $2$-variable weighted shift $\mathbf{T\equiv (}T_{1},T_{2})$
is horizontally flat (resp. vertically flat) if $\alpha
_{(k_{1},k_{2})}=\alpha _{(1,1)}$ for all $k_{1},k_{2}\geq 1$ (resp. $\beta
_{(k_{1},k_{2})}=\beta _{(1,1)}$ for all $k_{1},k_{2}\geq 1$). \ We say that 
$\mathbf{T}$ is flat if $\mathbf{T}$ is horizontally and vertically flat
(cf. Figure \ref{Figure 1}), and we say that $\mathbf{T}$ is symmetrically
flat if $\mathbf{T}$ is flat and $\alpha _{11}=\beta _{11}$. %
\setlength{\unitlength}{1mm} \psset{unit=1mm} 
\begin{figure}[th]
\begin{center}
\begin{picture}(140,138)

\psline{->}(20,20)(135,20)
\psline(20,40)(125,40)
\psline(20,60)(125,60)
\psline(20,80)(125,80)
\psline(20,100)(125,100)
\psline(20,120)(125,120)
\psline{->}(20,20)(20,135)
\psline(40,20)(40,125)
\psline(60,20)(60,125)
\psline(80,20)(80,125)
\psline(100,20)(100,125)
\psline(120,20)(120,125)

\put(11,16){\footnotesize{$(0,0)$}}
\put(36,16){\footnotesize{$(1,0)$}}
\put(57,16){\footnotesize{$\cdots$}}
\put(75,16){\footnotesize{$(k_{1},0)$}}
\put(93,16){\footnotesize{$(k_{1}+1,0)$}}
\put(117,16){\footnotesize{$\cdots$}}

\put(27,21){\footnotesize{$\alpha_{0,0}$}}
\put(47,21){\footnotesize{$\alpha_{1,0}$}}
\put(67,21){\footnotesize{$\cdots$}}
\put(87,21){\footnotesize{$\alpha_{k_{1},0}$}}
\put(105,21){\footnotesize{$\alpha_{k_{1}+1,0}$}}
\put(124,21){\footnotesize{$\cdots$}}

\put(27,41){\footnotesize{$\alpha_{0,1}$}}
\put(47,41){\footnotesize{$1$}}
\put(67,41){\footnotesize{$\cdots$}}
\put(87,41){\footnotesize{$1$}}
\put(107,41){\footnotesize{$1$}}
\put(124,41){\footnotesize{$\cdots$}}

\put(27,61){\footnotesize{$\alpha_{0,2}$}}
\put(47,61){\footnotesize{$1$}}
\put(67,61){\footnotesize{$\cdots$}}
\put(87,61){\footnotesize{$1$}}
\put(107,61){\footnotesize{$1$}}
\put(124,61){\footnotesize{$\cdots$}}

\put(27,81){\footnotesize{$\cdots$}}
\put(47,81){\footnotesize{$\cdots$}}
\put(67,81){\footnotesize{$\cdots$}}
\put(87,81){\footnotesize{$\cdots$}}
\put(107,81){\footnotesize{$\cdots$}}
\put(124,81){\footnotesize{$\cdots$}}

\put(27,101){\footnotesize{$\alpha_{0,k_{2}}$}}
\put(47,101){\footnotesize{$1$}}
\put(67,101){\footnotesize{$\cdots$}}
\put(87,101){\footnotesize{$1$}}
\put(107,101){\footnotesize{$1$}}
\put(124,101){\footnotesize{$\cdots$}}

\put(27,121){\footnotesize{$\alpha_{0,k_{2}+1}$}}
\put(47,121){\footnotesize{$1$}}
\put(67,121){\footnotesize{$\cdots$}}
\put(87,121){\footnotesize{$1$}}
\put(107,121){\footnotesize{$1$}}
\put(124,121){\footnotesize{$\cdots$}}

\psline{->}(70,10)(90,10)
\put(79,6){$\rm{T}_1$}

\put(11,38){\footnotesize{$(0,1)$}}
\put(11,58){\footnotesize{$(0,2)$}}
\put(14,78){\footnotesize{$\vdots$}}
\put(10,98){\footnotesize{$(0,k_{2})$}}
\put(5,118){\footnotesize{$(0,k_{2}+1)$}}

\psline{->}(10, 70)(10,90)
\put(5,80){$\rm{T}_2$}

\put(20,28){\footnotesize{$\beta_{0,0}$}}
\put(20,48){\footnotesize{$\beta_{0,1}$}}
\put(20,68){\footnotesize{$\beta_{0,2}$}}
\put(22,88){\footnotesize{$\vdots$}}
\put(20,108){\footnotesize{$\beta_{0,k_{2}}$}}
\put(22,128){\footnotesize{$\vdots$}}

\put(40,28){\footnotesize{$\beta_{1,0}$}}
\put(42,48){\footnotesize{$b$}}
\put(42,68){\footnotesize{$b$}}
\put(42,88){\footnotesize{$\vdots$}}
\put(42,108){\footnotesize{$b$}}
\put(42,128){\footnotesize{$\vdots$}}

\put(62,28){\footnotesize{$\vdots$}}
\put(62,48){\footnotesize{$\vdots$}}
\put(62,68){\footnotesize{$\vdots$}}
\put(62,88){\footnotesize{$\vdots$}}
\put(62,108){\footnotesize{$\vdots$}}
\put(62,128){\footnotesize{$\vdots$}}

\put(80,28){\footnotesize{$\beta_{k_{1},0}$}}
\put(82,48){\footnotesize{$b$}}
\put(82,68){\footnotesize{$b$}}
\put(82,88){\footnotesize{$\vdots$}}
\put(82,108){\footnotesize{$b$}}
\put(82,128){\footnotesize{$\vdots$}}

\put(100,28){\footnotesize{$\beta_{k_{1}+1,0}$}}
\put(102,48){\footnotesize{$b$}}
\put(102,68){\footnotesize{$b$}}
\put(102,88){\footnotesize{$\vdots$}}
\put(102,108){\footnotesize{$b$}}
\put(102,128){\footnotesize{$\vdots$}}

\put(122,28){\footnotesize{$\vdots$}}
\put(122,48){\footnotesize{$\vdots$}}
\put(122,68){\footnotesize{$\vdots$}}
\put(122,88){\footnotesize{$\vdots$}}
\put(122,108){\footnotesize{$\vdots$}}
\put(122,128){\footnotesize{$\vdots$}}

\put(49,40){\pscircle*(0,0){.75}}
\put(89,40){\pscircle*(0,0){.75}}
\put(109,40){\pscircle*(0,0){.75}}

\put(49,60){\pscircle*(0,0){.75}}
\put(89,60){\pscircle*(0,0){.75}}
\put(109,60){\pscircle*(0,0){.75}}

\put(49,100){\pscircle*(0,0){.75}}
\put(89,100){\pscircle*(0,0){.75}}
\put(109,100){\pscircle*(0,0){.75}}

\put(49,120){\pscircle*(0,0){.75}}
\put(89,120){\pscircle*(0,0){.75}}
\put(109,120){\pscircle*(0,0){.75}}

\psdots*[dotstyle=triangle*,dotangle=270](40.5,49)(40.5,69)(40.5,109)

\psdots*[dotstyle=triangle*,dotangle=270](80.5,49)(80.5,69)(80.5,109)

\psdots*[dotstyle=triangle*,dotangle=270](100.5,49)(100.5,69)(100.5,109)

\end{picture}
\end{center}
\caption{Weight diagram of a flat $2$-variable weighted shift (with round
dots for horizontal flatness, triangular dots for vertical flatness)}
\label{Figure 1}
\end{figure}
\end{definition}

\begin{lemma}
(\cite{bridge})\label{joint hypo}(Six-point Test) \ Let $\mathbf{T\equiv (}%
T_{1},T_{2})$ be a $2$-variable weighted shift, with weight sequences $%
\alpha $ and $\beta $. \ Then 
\begin{equation*}
\begin{tabular}{l}
$\lbrack \mathbf{T}^{\ast },\mathbf{T]\geq }0\Leftrightarrow (([T_{j}^{\ast
},T_{i}]e_{\mathbf{k+\varepsilon }_{j}},e_{\mathbf{k+\varepsilon }%
_{i}}))_{i,j=1}^{2}\geq 0\text{ (all }\mathbf{k}\in \mathbf{Z}_{+}^{2}\text{)%
}$ \\ 
\\ 
$\Leftrightarrow \left( 
\begin{array}{cc}
\alpha _{\mathbf{k}+\mathbf{\varepsilon }_{1}}^{2}-\alpha _{\mathbf{k}}^{2}
& \alpha _{\mathbf{k}+\mathbf{\varepsilon }_{2}}\beta _{\mathbf{k}+\mathbf{%
\varepsilon }_{1}}-\alpha _{\mathbf{k}}\beta _{\mathbf{k}} \\ 
\alpha _{\mathbf{k}+\mathbf{\varepsilon }_{2}}\beta _{\mathbf{k}+\mathbf{%
\varepsilon }_{1}}-\alpha _{\mathbf{k}}\beta _{\mathbf{k}} & \beta _{\mathbf{%
k}+\mathbf{\varepsilon }_{2}}^{2}-\beta _{\mathbf{k}}^{2}%
\end{array}%
\right) \geq 0\text{ (all }\mathbf{k}\in \mathbf{Z}_{+}^{2}\text{).}$%
\end{tabular}%
\end{equation*}
\end{lemma}

\setlength{\unitlength}{1mm} \psset{unit=1mm}

\begin{figure}[th]
\begin{center}
\begin{picture}(90,60)

\psline(20,10)(60,10)
\psline(20,30)(40,30)
\psline(20,10)(20,50)
\psline(40,10)(40,30)

\put(20,10){\pscircle*(0,0){1}}
\put(40,10){\pscircle*(0,0){1}}
\put(60,10){\pscircle*(0,0){1}}
\put(20,30){\pscircle*(0,0){1}}
\put(40,30){\pscircle*(0,0){1}}
\put(20,50){\pscircle*(0,0){1}}

\put(13,6){\footnotesize{$(k_1,k_2)$}}
\put(31,6){\footnotesize{$(k_1+1,k_2)$}}
\put(54,6){\footnotesize{$(k_1+2,k_2)$}}

\put(27,12){\footnotesize{$\alpha_{k_1,k_2}$}}
\put(47,12){\footnotesize{$\alpha_{k_1+1,k_2}$}}

\put(27,32){\footnotesize{$\alpha_{k_1,k_2+1}$}}

\put(3,29){\footnotesize{$(k_1,k_2+1)$}}
\put(3,49){\footnotesize{$(k_1, k_2+2)$}}
\put(41,29){\footnotesize{$(k_1+1, k_2+1)$}}

\put(20,19){\footnotesize{$\beta_{k_1,k_2}$}}
\put(20,39){\footnotesize{$\beta_{k_1,k_2+1}$}}

\put(40,19){\footnotesize{$\beta_{k_1+1,k_2}$}}

\end{picture}
\end{center}
\caption{{}Weight diagram for the Six-point Test}
\label{Figure 2}
\end{figure}
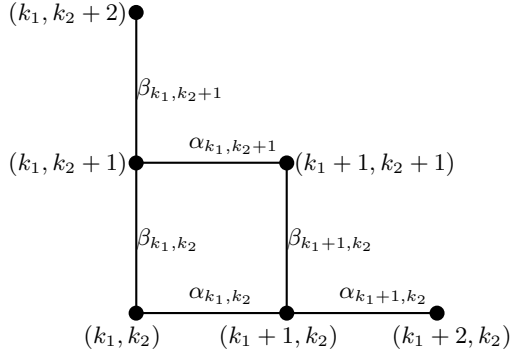

\begin{theorem}
\label{thm33}Let $\mathbf{T\equiv (}T_{1},T_{2})$ be a commuting, hyponormal 
$2$-variable weighted shift.\newline
(i) \ If $T_{1}$ is quadratically hyponormal and $\alpha _{(k_{1},k_{2})+%
\mathbf{\varepsilon }_{1}}=\alpha _{(k_{1},k_{2})}$ for some $%
k_{1},k_{2}\geq 1$, then $(T_{1},T_{2}(I\otimes U_{+}^{k_{2}-1}))$ is
horizontally flat.\newline
(ii) \ If, instead, $T_{2}$ is quadratically hyponormal and $\beta
_{(k_{1},k_{2})+\mathbf{\varepsilon }_{2}}=\beta _{(k_{1},k_{2})}$ for some $%
k_{1},k_{2}\geq 1$, then $(T_{1}(U_{+}^{k_{1}-1}\otimes I),T_{2})$ is
vertically flat.
\end{theorem}

\begin{proof}
Without loss of generality, we only prove (i). \ Consider the restricted
weight diagram based at$\ (k_{1},k_{2})$ (see Figure \ref{propagation}).%
\newline
\setlength{\unitlength}{1mm} \psset{unit=1mm}

\begin{figure}[th]
\begin{center}
\begin{picture}(90,50)

\psline(20,10)(82,10)
\psline(20,30)(82,30)
\psline(20,10)(20,50)
\psline(40,10)(40,40)
\psline(60,10)(60,40)
\psline(80,10)(80,40)

\psline{->}(36,3)(56,3)
\put(45,-0.5){$\rm{T}_1$}

\psline{->}(0, 19)(0,39)
\put(-6,30){$\rm{T}_2$}

\put(30,10){\pscircle*(0,0){1}}
\put(50,10){\pscircle*(0,0){1}}

\put(13,6){\footnotesize{$(k_1,k_2)$}}
\put(31,6){\footnotesize{$(k_1+1,k_2)$}}
\put(54,6){\footnotesize{$(k_1+2,k_2)$}}
\put(74,6){\footnotesize{$(k_1+3,k_2)$}}

\put(27,12){\footnotesize{$\alpha_{k_1,k_2}$}}
\put(45,12){\footnotesize{$\alpha_{k_1+1,k_2}$}}
\put(65,12){\footnotesize{$\alpha_{k_1+2,k_2}$}}

\put(26,32){\footnotesize{$\alpha_{k_1,k_2+1}$}}
\put(44,32){\footnotesize{$\alpha_{k_1+1,k_2+1}$}}

\put(4,29){\footnotesize{$(k_1,k_2+1)$}}
\put(4,49){\footnotesize{$(k_1,k_2+2)$}}

\put(20,19){\footnotesize{$\beta_{k_1,k_2}$}}
\put(20,39){\footnotesize{$\beta_{k_1,k_2+1}$}}

\put(40,19){\footnotesize{$\beta_{k_1+1,k_2}$}}

\end{picture}
\end{center}
\caption{Weight diagram of the $2$-variable weighted shift in Theorem \ref%
{thm33} (the two solid black dots represent equal weights)}
\label{propagation}
\end{figure}
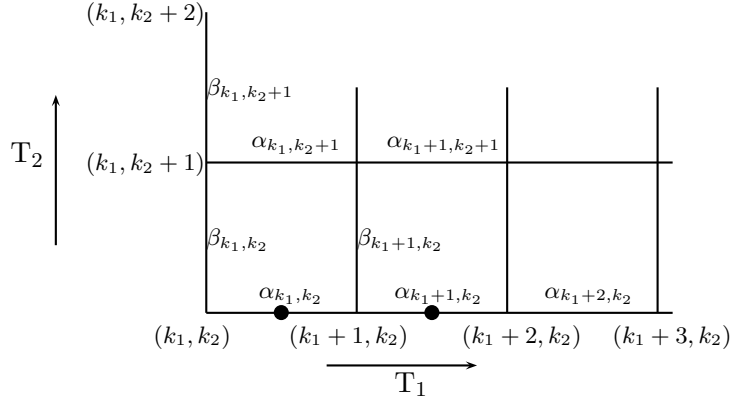
Recall that, by joint hyponormality, we have%
\begin{equation*}
\left( 
\begin{array}{cc}
\alpha _{(k_{1},k_{2})+\mathbf{\varepsilon }_{1}}^{2}-\alpha
_{(k_{1},k_{2})}^{2} & \alpha _{(k_{1},k_{2})+\mathbf{\varepsilon }%
_{2}}\beta _{(k_{1},k_{2})+\mathbf{\varepsilon }_{1}}-\beta
_{(k_{1},k_{2})}\alpha _{(k_{1},k_{2})} \\ 
\alpha _{(k_{1},k_{2})+\mathbf{\varepsilon }_{2}}\beta _{(k_{1},k_{2})+%
\mathbf{\varepsilon }_{1}}-\beta _{(k_{1},k_{2})}\alpha _{(k_{1},k_{2})} & 
\beta _{(k_{1},k_{2})+\mathbf{\varepsilon }_{2}}^{2}-\beta
_{(k_{1},k_{2})}^{2}%
\end{array}%
\right) \geq 0.
\end{equation*}%
Since $\alpha _{(k_{1},k_{2})+\mathbf{\varepsilon }_{1}}=\alpha
_{(k_{1},k_{2})}$, \ it follows that%
\begin{equation}
\alpha _{(k_{1},k_{2})+\mathbf{\varepsilon }_{2}}\beta _{(k_{1},k_{2})+%
\mathbf{\varepsilon }_{1}}=\beta _{(k_{1},k_{2})}\alpha _{(k_{1},k_{2})}.
\label{1}
\end{equation}%
By the commuting property (\ref{commuting}),%
\begin{equation}
\alpha _{(k_{1},k_{2})}\beta _{(k_{1},k_{2})+\mathbf{\varepsilon }%
_{1}}=\alpha _{(k_{1},k_{2})+\mathbf{\varepsilon }_{2}}\beta
_{(k_{1},k_{2})}.  \label{2}
\end{equation}%
Therefore 
\begin{equation*}
\begin{tabular}{l}
$\alpha _{(k_{1},k_{2})+\mathbf{\varepsilon }_{2}}^{2}\beta
_{(k_{1},k_{2})}=\alpha _{(k_{1},k_{2})+\mathbf{\varepsilon }_{2}}(\alpha
_{(k_{1},k_{2})+\mathbf{\varepsilon }_{2}}\beta _{(k_{1},k_{2})})=\alpha
_{(k_{1},k_{2})+\mathbf{\varepsilon }_{2}}(\alpha _{(k_{1},k_{2})}\beta
_{(k_{1},k_{2})+\mathbf{\varepsilon }_{1}})$ (by (\ref{2})) \\ 
\\ 
$=\alpha _{(k_{1},k_{2})}(\alpha _{(k_{1},k_{2})+\mathbf{\varepsilon }%
_{2}}\beta _{(k_{1},k_{2})+\mathbf{\varepsilon }_{1}})=\alpha
_{(k_{1},k_{2})}(\beta _{(k_{1},k_{2})}\alpha _{(k_{1},k_{2})})$ (by (\ref{1}%
)).%
\end{tabular}%
\end{equation*}%
Thus, $\alpha _{(k_{1},k_{2})+\mathbf{\varepsilon }_{2}}^{2}\beta
_{(k_{1},k_{2})}=\alpha _{(k_{1},k_{2})}(\beta _{(k_{1},k_{2})}\alpha
_{(k_{1},k_{2})})$, which implies that $\alpha _{(k_{1},k_{2})+\mathbf{%
\varepsilon }_{2}}=\alpha _{(k_{1},k_{2})}$. \ We now recall Theorem \ref%
{quadhypo3}, which says that flatness can be propagated to the right, that
is, $\alpha _{(k_{1},k_{2})+\mathbf{\varepsilon }_{1}}=\alpha
_{(k_{1},k_{2})+2\mathbf{\varepsilon }_{1}}$. \ It follows that $\alpha
_{(k_{1},k_{2})+\mathbf{\varepsilon }_{1}+\mathbf{\varepsilon }_{2}}=\alpha
_{(k_{1},k_{2})+\mathbf{\varepsilon }_{2}}$, and then two equal weights
occurs at level $k_{2}+1$, which then implies $\alpha _{(k_{1},k_{2})+2%
\mathbf{\varepsilon }_{2}}=\alpha _{(k_{1},k_{2})+\mathbf{\varepsilon }%
_{2}}=\alpha _{(k_{1},k_{2})}$. \ It is now easy to see that for every level 
$\ell \geq k_{2}$ we must have $\alpha _{(k_{1},\ell )}=\alpha
_{(k_{1},k_{2})}$ $($all $k_{1}\geq 1)$. \ Using Theorem \ref{quadhypo3} to
propagate these equalities to the left, we eventually conclude that%
\begin{equation*}
\alpha _{(k_{1},\ell )}=\alpha _{(1,k_{2})}\;\;(k_{1}\geq 1,\ell \geq k_{2}).
\end{equation*}%
We thus obtain that $(T_{1},T_{2})|_{\bigvee \{e_{(k_{1},\ell )}:\text{ }%
k_{1}\geq 1,\text{ }\ell \geq k_{2}\}}$ is unitarily equivalent to $(\alpha
_{(1,k_{2})}U_{+}\otimes I,I\otimes W_{\eta })$, where $\eta _{k}:=\beta
_{1,k+k_{2}}(k\geq 0)$. \ This can be rephrased as saying that $%
(T_{1},T_{2}(I\otimes U_{+}^{k_{2}-1}))$ is horizontally flat, as desired.
\end{proof}

\begin{remark}
\label{rem34}The proof of Theorem \ref{thm33} shows that for $\mathbf{%
T\equiv (}T_{1},T_{2})$ commuting and hyponormal, and for $k_{1},k_{2}\geq 0$%
, 
\begin{equation}
\alpha _{(k_{1},k_{2})+\mathbf{\varepsilon }_{1}}=\alpha
_{(k_{1},k_{2})}\Rightarrow \beta _{(k_{1},k_{2})}=\beta _{(k_{1},k_{2})+%
\mathbf{\varepsilon }_{1}}  \label{rmk34}
\end{equation}%
(by (\ref{1}) and (\ref{2})). \ Moreover, if $k_{2}\geq 1$, 
\begin{equation*}
\alpha _{(k_{1},k_{2})+\mathbf{\varepsilon }_{1}}=\alpha _{(k_{1},k_{2})}%
\text{ and }\alpha _{(k_{1},k_{2})+\mathbf{\varepsilon }_{1}-\mathbf{%
\varepsilon }_{2}}=\alpha _{(k_{1},k_{2})-\mathbf{\varepsilon }%
_{2}}\Rightarrow \alpha _{(k_{1},k_{2})+\mathbf{\varepsilon }_{1}}=\alpha
_{(k_{1},k_{2})+\mathbf{\varepsilon }_{1}-\mathbf{\varepsilon }_{2}}.
\end{equation*}
\end{remark}

\begin{remark}
\label{rem35}The proof of Theorem \ref{thm33} also reveals that asking $%
\mathbf{T}\equiv (T_{1},T_{2})$ to be jointly hyponormal is significantly
stronger than asking both $T_{1}$ and $T_{2}$ to be hyponormal. \ For,
consider the $2$-variable weighted shift whose weight diagram is given by
Figure \ref{totallyflat}. \ \setlength{\unitlength}{1mm} \psset{unit=1mm} 
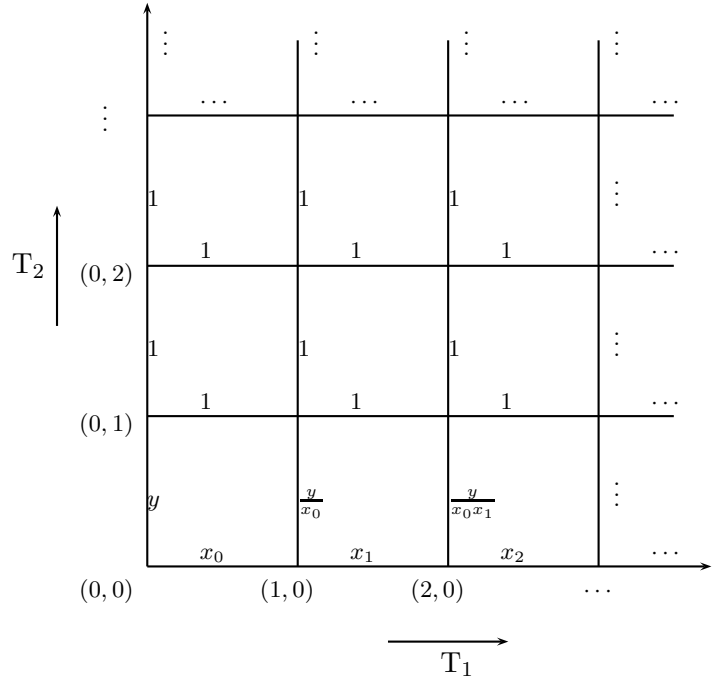
\begin{figure}[th]
\begin{center}
\begin{picture}(90,95)

\psline{->}(20,20)(95,20)
\psline(20,40)(90,40)
\psline(20,60)(90,60)
\psline(20,80)(90,80)
\psline{->}(20,20)(20,95)
\psline(40,20)(40,90)
\psline(60,20)(60,90)
\psline(80,20)(80,90)

\put(11,16){\footnotesize{$(0,0)$}}
\put(35,16){\footnotesize{$(1,0)$}}
\put(55,16){\footnotesize{$(2,0)$}}
\put(78,16){\footnotesize{$\cdots$}}

\put(27,21){\footnotesize{$x_0$}}
\put(47,21){\footnotesize{$x_1$}}
\put(67,21){\footnotesize{$x_2$}}
\put(87,21){\footnotesize{$\cdots$}}

\put(27,41){\footnotesize{$1$}}
\put(47,41){\footnotesize{$1$}}
\put(67,41){\footnotesize{$1$}}
\put(87,41){\footnotesize{$\cdots$}}

\put(27,61){\footnotesize{$1$}}
\put(47,61){\footnotesize{$1$}}
\put(67,61){\footnotesize{$1$}}
\put(87,61){\footnotesize{$\cdots$}}

\put(27,81){\footnotesize{$\cdots$}}
\put(47,81){\footnotesize{$\cdots$}}
\put(67,81){\footnotesize{$\cdots$}}
\put(87,81){\footnotesize{$\cdots$}}

\put(82,28){\footnotesize{$\vdots$}}
\put(82,48){\footnotesize{$\vdots$}}
\put(82,68){\footnotesize{$\vdots$}}
\put(82,88){\footnotesize{$\vdots$}}

\psline{->}(52,10)(68,10)
\put(59,6){$\rm{T}_1$}

\put(11,38){\footnotesize{$(0,1)$}}
\put(11,58){\footnotesize{$(0,2)$}}
\put(14,78){\footnotesize{$\vdots$}}

\psline{->}(8, 52)(8,68)
\put(2,59){$\rm{T}_2$}

\put(20,28){\footnotesize{$y$}}
\put(20,48){\footnotesize{$1$}}
\put(20,68){\footnotesize{$1$}}
\put(22,88){\footnotesize{$\vdots$}}

\put(40,28){\footnotesize{$\frac{y}{x_0}$}}
\put(40,48){\footnotesize{$1$}}
\put(40,68){\footnotesize{$1$}}
\put(42,88){\footnotesize{$\vdots$}}

\put(60,28){\footnotesize{$\frac{y}{x_0x_1}$}}
\put(60,48){\footnotesize{$1$}}
\put(60,68){\footnotesize{$1$}}
\put(62,88){\footnotesize{$\vdots$}}

\end{picture}
\end{center}
\caption{Weight diagram of the 2-variable weighted shift in Remark \ref%
{rem35}}
\label{totallyflat}
\end{figure}
In \cite[Theorem 5.2]{CuYo1}, we established that in the case when $%
\left\Vert W_{\alpha }\right\Vert \leq 1$, $\mathbf{T}$ is subnormal if and
only if $\mathbf{T}$ is hyponormal. \ Thus, a necessary condition for the
hyponormality of $\mathbf{T}$ is the subnormality of $W_{0}:=shift(\alpha
_{00},\alpha _{10},\cdots )$. \ For $0<a<1$, let $x_{0}\equiv x_{1}:=a$ and
let $x_{k}:=1\;(k\geq 2)$. \ Clearly $W_{0}$ is hyponormal and not
subnormal, and if we take $0<y\leq a^{2}$ we can guarantee that each of $%
T_{1}$ and $T_{2}$ is hyponormal, yet $\mathbf{T}$ is not. \ An alternative
way to see this is observe that if $\mathbf{T}$ were hyponormal then $\alpha
_{01}$ would equal $a$, since $\alpha _{00}=\alpha _{10}$.
\end{remark}

We will now show that Theorem \ref{thm33} is optimal in the following sense:
the propagation does not necessarily extend either to the left ($0$-th
column) or down (below $k_{2}$-th level). \ To demonstrate this optimality,
we first introduce the class of Bergman-like weighted shifts.

\begin{definition}
For $\ell \geq 1$, the Bergman-like weighted shift on $\ell ^{2}(\mathbb{Z}%
_{+})$ is $B_{+}^{(\ell )}:=$ $shift(\{\sqrt{\ell -\frac{1}{k+2}}:k\geq 0\})$%
; that is, 
\begin{equation*}
B_{+}^{(\ell )}e_{k}:=\sqrt{\ell -\frac{1}{k+2}}e_{k+1}\;\;(k\geq 0).
\end{equation*}%
In particular, $B_{+}^{(1)}\equiv B_{+}:=shift(\sqrt{\frac{1}{2}},\sqrt{%
\frac{2}{3}},\sqrt{\frac{3}{4}},\cdots )$ is the Bergman shift.
\end{definition}

\begin{remark}
(i) $\ B_{+}$ is subnormal with Berger measure $d\xi (s):=ds$ on $[0,1]$.%
\newline
(ii) \ (\cite{CuYo3}) \ $B_{+}^{(2)}$ is subnormal with Berger measure $d\xi
(s):=\frac{sds}{\pi \sqrt{2s-s^{2}}}$ on $[0,2]$.
\end{remark}

Our next step is to show that $B_{+}^{(\ell )}$ $(\ell \geq 1)$ is always $2$%
-hyponormal. \ To this end, we need two preliminary results.

\begin{lemma}
\label{Nested Determinants Test}(Nested Determinants Test; Special Case) \
If $a>0$\ and $\det \left( 
\begin{array}{cc}
a & b \\ 
b & c%
\end{array}%
\right) >0$, then%
\begin{equation*}
\left( 
\begin{array}{ccc}
a & b & c \\ 
b & c & d \\ 
c & d & e%
\end{array}%
\right) \geq 0\Leftrightarrow \det \left( 
\begin{array}{ccc}
a & b & c \\ 
b & c & d \\ 
c & d & e%
\end{array}%
\right) \geq 0.
\end{equation*}
\end{lemma}

\begin{proof}
Straightforward from Choleski's Algorithm (\cite{Atk}).
\end{proof}

\begin{lemma}
(\cite{QHWS})\label{thm:k-hyponormal}Let $W_{\alpha }e_{k}=\alpha
_{k}e_{k+1} $ $(k\geq 0)$ be a hyponormal weighted shift. \ The following
statements are equivalent:\newline
(i) $\ W_{\alpha }$ is $2$-hyponormal.\newline
(ii) \ The matrix 
\begin{equation*}
(([W_{\alpha }^{\ast j},W_{\alpha }^{i}]e_{k+j},e_{k+i}))_{i,j=1}^{2}
\end{equation*}%
is positive semi-definite for all $k\geq -1.$\newline
(iii) \ The matrix 
\begin{equation*}
(\gamma _{k}\gamma _{k+i+j}-\gamma _{k+i}\gamma _{k+j})_{i,j=1}^{2}
\end{equation*}%
is positive semi-definite for all $k\geq 0,$ where as usual $\gamma _{0}=1,$ 
$\gamma _{n}=\alpha _{0}^{2}\cdot \cdots \cdot \alpha _{n-1}^{2}$ $(n\geq 1)$%
.\newline
(iv) \ The Hankel matrix 
\begin{equation*}
H(2;k):=(\gamma _{k+i+j-2})_{i,j=1}^{3}
\end{equation*}%
is positive semi-definite for all $k\geq 0$.
\end{lemma}

We now use symbolic manipulation to prove the following result.

\begin{theorem}
\label{2hypo}All Bergman-like weighted shifts $B_{+}^{(\ell )}$ $($all $\ell
\geq 1)$ are $2$-hyponormal.
\end{theorem}

\begin{proof}
By Lemma \ref{Nested Determinants Test} and Lemma \ref{thm:k-hyponormal}, to
prove that $B_{+}^{(\ell )}$ is $2$-hyponormal it suffices to see that $\det
H(2;k)>0$ for all $k\geq 0$. \ Now, 
\begin{eqnarray*}
\det H(2;k) &=&\gamma _{k}^{3}\det \left( 
\begin{array}{ccc}
1 & \alpha _{k}^{2} & \alpha _{k}^{2}\alpha _{k+1}^{2} \\ 
&  &  \\ 
\alpha _{k}^{2} & \alpha _{k}^{2}\alpha _{k+1}^{2} & \alpha _{k}^{2}\alpha
_{k+1}^{2}\alpha _{k+2}^{2} \\ 
&  &  \\ 
\alpha _{k}^{2}\alpha _{k+1}^{2} & \alpha _{k}^{2}\alpha _{k+1}^{2}\alpha
_{k+2}^{2} & \alpha _{k}^{2}\alpha _{k+1}^{2}\alpha _{k+2}^{2}\alpha
_{k+3}^{2}%
\end{array}%
\right) \\
&& \\
&=&\gamma _{k}^{3}\frac{2(\ell +1)((k+2)\ell -1)^{2}((k+3)\ell -1)}{%
(k+2)^{3}(k+3)^{3}(k+4)^{2}(k+5)}>0,
\end{eqnarray*}%
as desired.
\end{proof}

\begin{corollary}
For every $\ell \geq 1$, the Bergman-like weighted shift $B_{+}^{(\ell )}$
is quadratically hyponormal.

\begin{remark}
In (\cite{CPY}), we prove a much stronger result: all Bergman-like weighted
shifts $B_{+}^{(\ell )}$ $($all $\ell \geq 1)$ are subnormal.
\end{remark}
\end{corollary}

Theorem \ref{Propagation of Induction} below says that the amount of
propagation provided by Theorem \ref{thm33} is maximum; briefly, we say that
Theorem \ref{thm33} is optimal. \ Observe that for the $2$-variable weighted
shift in Figure \ref{Figure 5}, we have $\alpha _{(k_{1},k_{2})+\mathbf{%
\varepsilon }_{1}}=\alpha _{(k_{1},k_{2})}\;\;($all $k_{1}\geq 1,k_{2}\geq
2) $, yet $\alpha _{(k_{1},k_{2})}<\alpha _{(k_{1},k_{2})+\mathbf{%
\varepsilon }_{1}}$ for all $k_{1}\geq 0$ and $k_{2}=0,1$ and $\alpha
_{(0,k_{2})}<\alpha _{(1,k_{2})}$ for all $k_{2}\geq 0$. \ In other words,
the trivial weight structure present in the subspace $\bigvee
\{e_{(k_{1},k_{2})}:k_{1}\geq 1,k_{2}\geq 2\}$ cannot be expanded either to
the left ($0$th column) or down (first row). \ First, we need an auxiliary
result, of independent interest.

\begin{lemma}
\label{lemofhypoflat}Consider the $2$-variable weighted shift $\mathbf{T}%
\equiv (T_{1},T_{2})$ given by Figure \ref{lem of induction}, where \newline
$shift(x_{0},x_{1},x_{2},\cdots )$ and $shift(y_{0},y_{1},y_{2},\cdots )$
are Bergman-like weighted shifts. \ Assume that \newline
$(T_{1},T_{2})\mathbf{|}_{\mathcal{M}}$ is jointly hyponormal, where $%
\mathcal{M}$ is the subspace associated to indices $\mathbf{k}$ with $%
k_{2}\geq 1$. \ Then there exists a Bergman-like weighted shift $%
shift(z_{0},z_{1},z_{2},\cdots )$ and a hyponormal weighted shift $W_{\beta
}:=shift(\beta _{0},\beta _{1},\beta _{2},\cdots )\;$\ $(\beta _{n}<\beta
_{n+1}$ for all $n\geq 0)$ such that $\mathbf{T}$ is jointly hyponormal.
\end{lemma}

\setlength{\unitlength}{1mm} \psset{unit=1mm} 
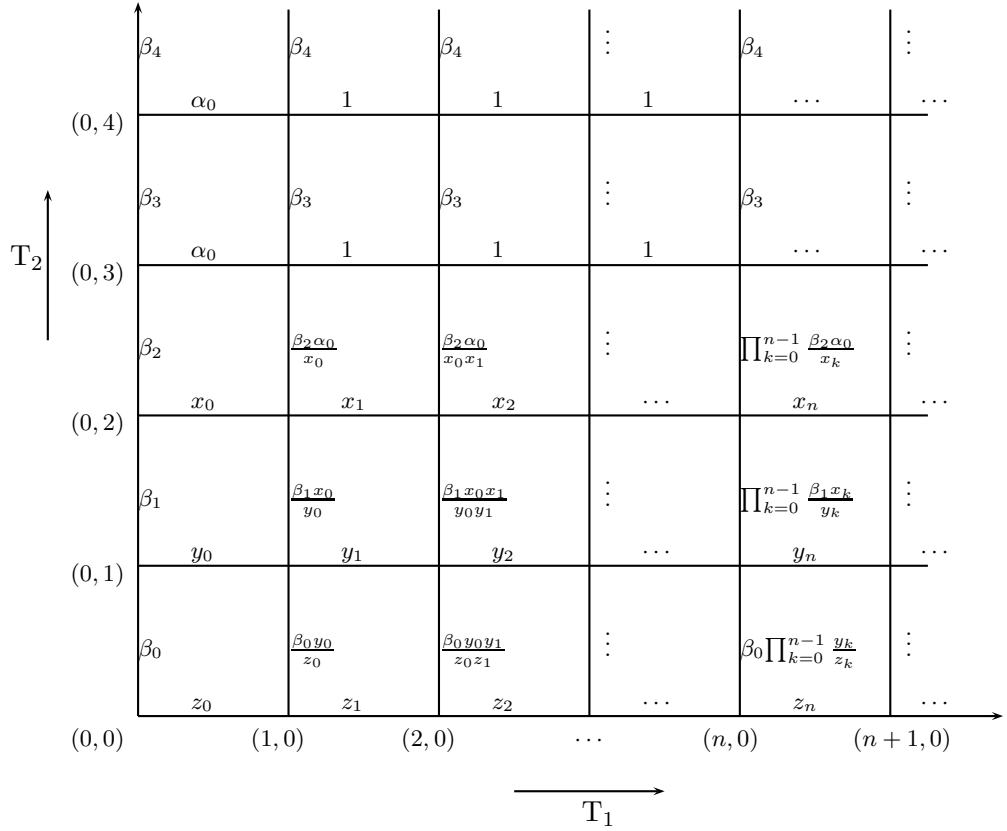
\begin{figure}[t]
\begin{center}
\begin{picture}(140,120)

\psline{->}(20,20)(135,20)
\psline(20,40)(125,40)
\psline(20,60)(125,60)
\psline(20,80)(125,80)
\psline(20,100)(125,100)
\psline{->}(20,20)(20,115)
\psline(40,20)(40,114)
\psline(60,20)(60,114)
\psline(80,20)(80,114)
\psline(100,20)(100,114)
\psline(120,20)(120,114)

\put(11,16){\footnotesize{$(0,0)$}}
\put(35,16){\footnotesize{$(1,0)$}}
\put(55,16){\footnotesize{$(2,0)$}}
\put(78,16){\footnotesize{$\cdots$}}
\put(95,16){\footnotesize{$(n,0)$}}
\put(115,16){\footnotesize{$(n+1,0)$}}

\put(27,21){\footnotesize{$z_{0}$}}
\put(47,21){\footnotesize{$z_{1}$}}
\put(67,21){\footnotesize{$z_{2}$}}
\put(87,21){\footnotesize{$\cdots$}}
\put(107,21){\footnotesize{$z_{n}$}}
\put(124,21){\footnotesize{$\cdots$}}

\put(27,41){\footnotesize{$y_{0}$}}
\put(47,41){\footnotesize{$y_{1}$}}
\put(67,41){\footnotesize{$y_{2}$}}
\put(87,41){\footnotesize{$\cdots$}}
\put(107,41){\footnotesize{$y_{n}$}}
\put(124,41){\footnotesize{$\cdots$}}

\put(27,61){\footnotesize{$x_{0}$}}
\put(47,61){\footnotesize{$x_{1}$}}
\put(67,61){\footnotesize{$x_{2}$}}
\put(87,61){\footnotesize{$\cdots$}}
\put(107,61){\footnotesize{$x_{n}$}}
\put(124,61){\footnotesize{$\cdots$}}

\put(27,81){\footnotesize{$\alpha_{0}$}}
\put(47,81){\footnotesize{$1$}}
\put(67,81){\footnotesize{$1$}}
\put(87,81){\footnotesize{$1$}}
\put(107,81){\footnotesize{$\cdots$}}
\put(124,81){\footnotesize{$\cdots$}}

\put(27,101){\footnotesize{$\alpha_{0}$}}
\put(47,101){\footnotesize{$1$}}
\put(67,101){\footnotesize{$1$}}
\put(87,101){\footnotesize{$1$}}
\put(107,101){\footnotesize{$\cdots$}}
\put(124,101){\footnotesize{$\cdots$}}

\psline{->}(70,10)(90,10)
\put(79,6){$\rm{T}_1$}

\put(11,38){\footnotesize{$(0,1)$}}
\put(11,58){\footnotesize{$(0,2)$}}
\put(11,78){\footnotesize{$(0,3)$}}
\put(11,98){\footnotesize{$(0,4)$}}

\psline{->}(8, 70)(8,90)
\put(3,80){$\rm{T}_2$}

\put(20,28){\footnotesize{$\beta_{0}$}}
\put(20,48){\footnotesize{$\beta_{1}$}}
\put(20,68){\footnotesize{$\beta_{2}$}}
\put(20,88){\footnotesize{$\beta_{3}$}}
\put(20,108){\footnotesize{$\beta_{4}$}}

\put(40,28){\footnotesize{$\frac{\beta_{0}y_{0}}{z_{0}}$}}
\put(40,48){\footnotesize{$\frac{\beta_{1}x_{0}}{y_{0}}$}}
\put(40,68){\footnotesize{$\frac{\beta_{2}\alpha_{0}}{x_{0}}$}}
\put(40,88){\footnotesize{$\beta_{3}$}}
\put(40,108){\footnotesize{$\beta_{4}$}}

\put(60,28){\footnotesize{$\frac{\beta_{0}y_{0}y_{1}}{z_{0}z_{1}}$}}
\put(60,48){\footnotesize{$\frac{\beta_{1}x_{0}x_{1}}{y_{0}y_{1}}$}}
\put(60,68){\footnotesize{$\frac{\beta_{2}\alpha_{0}}{x_{0}x_{1}}$}}
\put(60,88){\footnotesize{$\beta_{3}$}}
\put(60,108){\footnotesize{$\beta_{4}$}}

\put(82,28){\footnotesize{$\vdots$}}
\put(82,48){\footnotesize{$\vdots$}}
\put(82,68){\footnotesize{$\vdots$}}
\put(82,88){\footnotesize{$\vdots$}}
\put(82,108){\footnotesize{$\vdots$}}

\put(100,28){\footnotesize{$\beta_{0}{\prod_{k=0}^{n-1}\frac{y_{k}}{z_{k}}}$}}
\put(100,48){\footnotesize{${\prod_{k=0}^{n-1}\frac{\beta_{1}x_{k}}{y_{k}}}$}}
\put(100,68){\footnotesize{${\prod_{k=0}^{n-1}\frac{\beta_{2}\alpha_{0}}{x_{k}}}$}}
\put(100,88){\footnotesize{$\beta_{3}$}}
\put(100,108){\footnotesize{$\beta_{4}$}}

\put(122,28){\footnotesize{$\vdots$}}
\put(122,48){\footnotesize{$\vdots$}}
\put(122,68){\footnotesize{$\vdots$}}
\put(122,88){\footnotesize{$\vdots$}}
\put(122,108){\footnotesize{$\vdots$}}

\end{picture}
\end{center}
\caption{Weight diagram of the 2-variable weighted shift in Lemma \ref%
{lemofhypoflat}}
\label{lem of induction}
\end{figure}

\begin{proof}
Let%
\begin{equation*}
\begin{tabular}{l}
$shift(x_{0},x_{1},x_{2},\cdots )\equiv shift(\{\sqrt{p-\frac{1}{n+2}}:n\geq
0\})$, \\ 
$shift(y_{0},y_{1},y_{2},\cdots )\equiv shift(\{\sqrt{q-\frac{1}{n+2}}:n\geq
0\})$ and \\ 
$shift(z_{0},z_{1},z_{2},\cdots )\equiv shift(\{\sqrt{r-\frac{1}{n+2}}:n\geq
0\})$, for some integers $p<q<r$.%
\end{tabular}%
\end{equation*}%
Since the restriction of $(T_{1},T_{2})$ to $\vee
\{e_{(k_{1},k_{2})}:k_{2}\geq 1\}$ is jointly hyponormal, it suffices to
apply the Six-point Test (Lemma \ref{joint hypo}) to $\mathbf{k}=(n,0)$,
with $n\geq 0$.\newline
\textbf{Case 1}: $\ \mathbf{k}=(0,0)$. \ Here%
\begin{equation*}
\begin{tabular}{l}
$M(0,0):=\left( 
\begin{array}{cc}
z_{1}^{2}-z_{0}^{2} & \frac{y_{0}^{2}\beta _{0}}{z_{0}}-\beta _{0}\cdot z_{0}
\\ 
\frac{y_{0}^{2}\beta _{0}}{z_{0}}-\beta _{0}\cdot z_{0} & \beta
_{1}^{2}-\beta _{0}^{2}%
\end{array}%
\right) \geq 0$ \\ 
\\ 
$\Leftrightarrow z_{0}^{2}(z_{1}^{2}-z_{0}^{2})(\beta _{1}^{2}-\beta
_{0}^{2})\geq \beta _{0}^{2}(z_{0}^{2}-y_{0}^{2})^{2}$ \\ 
\\ 
$\Leftrightarrow (r-\frac{1}{2})(\beta _{1}^{2}-\beta _{0}^{2})\geq 6\beta
_{0}^{2}(r-q)^{2}$.%
\end{tabular}%
\end{equation*}%
If we choose $\beta _{0}$ such that $\beta _{0}\leq \beta _{1}$ and $\beta
_{1}^{2}\geq \frac{12(r-q)^{2}+2r-1}{2r-1}\beta _{0}^{2}$, we obtain $%
M(0,0)\geq 0$.\newline
\textbf{Case 2}: \ $\mathbf{k}=(n,0)$ $\;(n\geq 1)$. \ Here%
\begin{equation*}
\begin{tabular}{l}
$M(n,0):=\left( 
\begin{array}{cc}
z_{n+1}^{2}-z_{n}^{2} & y_{n}\beta _{0}\prod_{k=0}^{n}\frac{y_{k}}{z_{k}}%
-z_{n}\beta _{0}\prod_{k=0}^{n-1}\frac{y_{k}}{z_{k}} \\ 
y_{n}\beta _{0}\prod_{k=0}^{n}\frac{y_{k}}{z_{k}}-z_{n}\beta
_{0}\prod_{k=0}^{n-1}\frac{y_{k}}{z_{k}} & \beta _{1}^{2}\prod_{k=0}^{n-1}(%
\frac{x_{k}}{y_{k}})^{2}-\beta _{0}^{2}\prod_{k=0}^{n-1}(\frac{y_{k}}{z_{k}}%
)^{2}%
\end{array}%
\right) \geq 0$ \\ 
\\ 
$\Leftrightarrow (z_{n+1}^{2}-z_{n}^{2})(\beta _{1}^{2}\prod_{k=0}^{n-1}(%
\frac{x_{k}}{y_{k}})^{2}-\beta _{0}^{2}\prod_{k=0}^{n-1}(\frac{y_{k}}{z_{k}}%
)^{2})\geq \beta _{0}^{2}(y_{n}\prod_{k=0}^{n}\frac{y_{k}}{z_{k}}%
-z_{n}\prod_{k=0}^{n-1}\frac{y_{k}}{z_{k}})^{2}$ \\ 
\\ 
$\Leftrightarrow z_{n}^{2}(z_{n+1}^{2}-z_{n}^{2})(\beta
_{1}^{2}\prod_{k=0}^{n-1}(\frac{x_{k}}{y_{k}})^{2}-\beta
_{0}^{2}\prod_{k=0}^{n-1}(\frac{y_{k}}{z_{k}})^{2})\geq \beta
_{0}^{2}\prod_{k=0}^{n-1}(\frac{y_{k}}{z_{k}})^{2}(y_{n}^{2}-z_{n}^{2})^{2}$
\\ 
\\ 
$\Leftrightarrow z_{n}^{2}(z_{n+1}^{2}-z_{n}^{2})(\beta
_{1}^{2}\prod_{k=0}^{n-1}(\frac{x_{k}z_{k}}{y_{k}^{2}})^{2}-\beta
_{0}^{2})\geq \beta _{0}^{2}(y_{n}^{2}-z_{n}^{2})^{2}$%
\end{tabular}%
\end{equation*}%
If we choose $p,$ $q$ and $r$ such that $\frac{x_{k}z_{k}}{y_{k}^{2}}\geq 3$%
, then $\frac{r(n+2)-1}{(n+2)^{2}(n+3)}((\frac{12(r-q)^{2}+2r-1}{2r-1}%
)\prod_{k=0}^{n-1}(\frac{x_{k}z_{k}}{y_{k}^{2}})^{2}-1)\geq (r-q)^{2}$ $\;($%
all $n\geq 1)$, which implies $M(n,0)\geq 0$ $\;($all $n\geq 1)$. \ 

By Cases 1 and 2, it follows that $(T_{1},T_{2})$ is jointly hyponormal.
\end{proof}

\begin{theorem}
\label{Propagation of Induction}For every $k_{2}\geq 1$ and $0<\alpha _{0}<1$
there exist \newline
(i) a family $\{B_{+}^{(\ell _{i})}\}_{i=0}^{k_{2}-1}$ of Bergman-like
weighted shifts, and \newline
(ii) a subnormal weighted shift $W_{\beta }:=shift(\beta _{0},\beta
_{1},\beta _{2},\cdots )$ (with $\beta _{n}<\beta _{n+1}$ for all $n\geq 0$),%
\newline
such that the commuting $2$-variable weighted shift $\mathbf{T}$ with a
weight diagram whose first $k_{2}$ rows are $B_{+}^{(\ell _{0})},\cdots
,B_{+}^{(\ell _{k_{2}-1})}$, whose remaining rows are $S_{\alpha _{0}}$, and
whose $0$-th column is given by $W_{\beta }$, is (jointly) hyponormal (see
Figure \ref{Figure 5} for the case $k_{2}=2$).
\end{theorem}

\begin{proof}
We divide the proof into three cases, according to the value of $k_{2}$.%
\newline
\textbf{Case 1}: $k_{2}=1$. \ For $p\geq 1$ let $\alpha _{m,0}\equiv x_{m}:=%
\sqrt{p-\frac{1}{m+2}}\;\;(m\geq 0)$. \ Since the restriction of $%
(T_{1},T_{2})$ to $\bigvee \{e_{(k_{1},k_{2})}:k_{2}\geq 1\}$ must be
unitarily equivalent to $(S_{\alpha _{0}}\otimes I,I\otimes shift(\beta
_{1},\beta _{2},\cdots ))$, to guarantee the hyponormality of $(T_{1},T_{2})$
it suffices to apply the Six-point Test (Lemma \ref{joint hypo}) to $\mathbf{%
k}=(m,0)$, with $m\geq 0$.
\end{proof}

\setlength{\unitlength}{1mm} \psset{unit=1mm} 
\begin{figure}[tbp]
\begin{center}
\begin{picture}(140,120)

\psline{->}(20,20)(135,20)
\psline(20,40)(125,40)
\psline(20,60)(125,60)
\psline(20,80)(125,80)
\psline(20,100)(125,100)
\psline{->}(20,20)(20,115)
\psline(40,20)(40,114)
\psline(60,20)(60,114)
\psline(80,20)(80,114)
\psline(100,20)(100,114)
\psline(120,20)(120,114)

\put(11,16){\footnotesize{$(0,0)$}}
\put(35,16){\footnotesize{$(1,0)$}}
\put(55,16){\footnotesize{$(2,0)$}}
\put(78,16){\footnotesize{$\cdots$}}
\put(95,16){\footnotesize{$(m,0)$}}
\put(115,16){\footnotesize{$(m+1,0)$}}

\put(27,21){\footnotesize{$y_{0}$}}
\put(47,21){\footnotesize{$y_{1}$}}
\put(67,21){\footnotesize{$y_{2}$}}
\put(87,21){\footnotesize{$\cdots$}}
\put(107,21){\footnotesize{$y_{m}$}}
\put(124,21){\footnotesize{$\cdots$}}

\put(27,41){\footnotesize{$x_{0}$}}
\put(47,41){\footnotesize{$x_{1}$}}
\put(67,41){\footnotesize{$x_{2}$}}
\put(87,41){\footnotesize{$\cdots$}}
\put(107,41){\footnotesize{$x_{m}$}}
\put(124,41){\footnotesize{$\cdots$}}

\put(27,61){\footnotesize{$\alpha_{0}$}}
\put(47,61){\footnotesize{$1$}}
\put(67,61){\footnotesize{$1$}}
\put(87,61){\footnotesize{$\cdots$}}
\put(107,61){\footnotesize{$1$}}
\put(124,61){\footnotesize{$\cdots$}}

\put(27,81){\footnotesize{$\cdots$}}
\put(47,81){\footnotesize{$\cdots$}}
\put(67,81){\footnotesize{$\cdots$}}
\put(87,81){\footnotesize{$\cdots$}}
\put(107,81){\footnotesize{$\cdots$}}
\put(124,81){\footnotesize{$\cdots$}}

\put(27,101){\footnotesize{$\alpha_{0}$}}
\put(47,101){\footnotesize{$1$}}
\put(67,101){\footnotesize{$1$}}
\put(87,101){\footnotesize{$\cdots$}}
\put(107,101){\footnotesize{$1$}}
\put(124,101){\footnotesize{$\cdots$}}

\psline{->}(70,10)(90,10)
\put(79,6){$\rm{T}_1$}

\put(11,38){\footnotesize{$(0,1)$}}
\put(11,58){\footnotesize{$(0,2)$}}
\put(14,78){\footnotesize{$\vdots$}}
\put(11,98){\footnotesize{$(0,n)$}}

\psline{->}(10, 70)(10,90)
\put(5,80){$\rm{T}_2$}

\put(20,28){\footnotesize{$\beta_{0}$}}
\put(20,48){\footnotesize{$\beta_{1}$}}
\put(20,68){\footnotesize{$\beta_{2}$}}
\put(22,88){\footnotesize{$\vdots$}}
\put(20,108){\footnotesize{$\beta_{n}$}}

\put(40,28){\footnotesize{$\frac{\beta_{0}x_{0}}{y_{0}}$}}
\put(40,48){\footnotesize{$\frac{\alpha_{0}\beta_{1}}{x_{0}}$}}
\put(40,68){\footnotesize{$\beta_{2}$}}
\put(42,88){\footnotesize{$\vdots$}}
\put(40,108){\footnotesize{$\beta_{n}$}}

\put(60,28){\footnotesize{$\frac{\beta_{0}x_{0}x_{1}}{y_{0}y_{1}}$}}
\put(60,48){\footnotesize{$\frac{\alpha_{0}\beta_{1}}{x_{0}x_{1}}$}}
\put(60,68){\footnotesize{$\beta_{2}$}}
\put(62,88){\footnotesize{$\vdots$}}
\put(60,108){\footnotesize{$\beta_{n}$}}

\put(82,28){\footnotesize{$\vdots$}}
\put(82,48){\footnotesize{$\vdots$}}
\put(82,68){\footnotesize{$\vdots$}}
\put(82,88){\footnotesize{$\vdots$}}
\put(82,108){\footnotesize{$\vdots$}}

\put(100,28){\footnotesize{$\beta_{0}{\prod_{k=0}^{m-1}\frac{x_{k}}{y_{k}}}$}}
\put(100,48){\footnotesize{${\prod_{k=0}^{m-1}\frac{\alpha_{0}\beta_{1}}{x_{k}}}$}}
\put(100,68){\footnotesize{$\beta_{2}$}}
\put(102,88){\footnotesize{$\vdots$}}
\put(100,108){\footnotesize{$\beta_{n}$}}

\put(122,28){\footnotesize{$\vdots$}}
\put(122,48){\footnotesize{$\vdots$}}
\put(122,68){\footnotesize{$\vdots$}}
\put(122,88){\footnotesize{$\vdots$}}
\put(122,108){\footnotesize{$\vdots$}}

\end{picture}
\end{center}
\caption{Weight diagram of the $2$-variable weighted shift in Theorem \ref%
{Propagation of Induction}}
\label{Figure 5}
\end{figure}

\textbf{Subcase 1}: $\ \mathbf{k}=(0,0)$. \ Here we have 
\begin{eqnarray}
M(0,0) &:&=\left( 
\begin{array}{cc}
x_{1}^{2}-x_{0}^{2} & \frac{\alpha _{0}^{2}\beta _{0}}{x_{0}}-\beta _{0}x_{0}
\\ 
\frac{\alpha _{0}^{2}\beta _{0}}{x_{0}}-\beta _{0}x_{0} & \beta
_{1}^{2}-\beta _{0}^{2}%
\end{array}%
\right)  \notag \\
&&  \notag \\
&=&\left( 
\begin{array}{cc}
\frac{1}{6} & \frac{\beta _{0}}{x_{0}}(\alpha _{0}^{2}-x_{0}^{2}) \\ 
\frac{\beta _{0}}{x_{0}}(\alpha _{0}^{2}-x_{0}^{2}) & \beta _{1}^{2}-\beta
_{0}^{2}%
\end{array}%
\right) \geq 0  \notag \\
&&  \notag \\
&\Leftrightarrow &6\beta _{0}^{2}(\alpha _{0}^{2}-x_{0}^{2})^{2}\leq (\beta
_{1}^{2}-\beta _{0}^{2})x_{0}^{2},  \label{beta0eq2}
\end{eqnarray}%
which imposes a condition on $x_{0}$ and $\beta _{0}$.\newline
\textbf{Subcase 2}: \ $\mathbf{k}=(m,0)$, with $m\geq 1$. \ Fix $m\geq 1$
and let $P_{m}:=\prod_{k=0}^{m-1}x_{k}$. \ We then see that 
\begin{eqnarray*}
M(m,0) &:&=\left( 
\begin{array}{cc}
x_{m+1}^{2}-x_{m}^{2} & \frac{\alpha _{0}\beta _{0}}{x_{m}P_{m}}-x_{m}\frac{%
\alpha _{0}\beta _{0}}{P_{m}} \\ 
\frac{\alpha _{0}\beta _{0}}{x_{m}P_{m}}-x_{m}\frac{\alpha _{0}\beta _{0}}{%
P_{m}} & \beta _{1}^{2}-\frac{\alpha _{0}^{2}\beta _{0}^{2}}{P_{m}^{2}}%
\end{array}%
\right) \geq 0 \\
&& \\
&\Leftrightarrow &x_{m}^{2}(x_{m+1}^{2}-x_{m}^{2})(\beta
_{1}^{2}P_{m}^{2}-\alpha _{0}^{2}\beta _{0}^{2})\geq \alpha _{0}^{2}\beta
_{0}^{2}(1-x_{m}^{2})^{2} \\
&&
\end{eqnarray*}%
\begin{eqnarray}
&\Leftrightarrow &x_{m}^{2}(\beta _{1}^{2}P_{m}^{2}-\alpha _{0}^{2}\beta
_{0}^{2})\geq (m+2)(m+3)\alpha _{0}^{2}\beta _{0}^{2}(1-x_{m}^{2})^{2} 
\notag \\
&&  \notag \\
&\Leftrightarrow &\frac{\beta _{1}^{2}P_{m}^{2}}{\alpha _{0}^{2}\beta
_{0}^{2}}\geq 1+(m+2)(m+3)(\frac{1}{x_{m}}-x_{m})^{2}.  \label{beta0eq}
\end{eqnarray}%
We now let $p=3$, so that $x_{k}^{2}\geq 2\;\;($all $k\geq 0)$ and therefore 
$P_{m}\geq 2^{m}\;\;($all $m\geq 1)$. \ Since $\lim_{m\rightarrow \infty }%
\frac{2^{m}}{(m+2)(m+3)}=\infty $, it is clear that we can find $\beta _{0}$
sufficiently small so that both (\ref{beta0eq2}) and (\ref{beta0eq}) hold.

\begin{proof}
From Subcases 1 and 2, it follows that $\mathbf{T}$ is jointly hyponormal.%
\newline
\textbf{Case 2}: $k_{2}=2$. \ Here we let $p:=3$, $q:=18$, $\beta
_{2}:=4\beta _{1}$ and $\beta _{1}:=\frac{1}{\alpha _{0}}$, so that $\alpha
_{m,1}\equiv x_{m}:=\sqrt{p-\frac{1}{m+2}}\equiv \sqrt{3-\frac{1}{m+2}}$ and 
$\alpha _{m,0}\equiv y_{m}:=\sqrt{q-\frac{1}{m+2}}\equiv \sqrt{18-\frac{1}{%
m+2}}\;\;(m\geq 0)$. \ Since the restriction of $(T_{1},T_{2})$ to $\bigvee
\{e_{(k_{1},k_{2})}:k_{2}\geq 2\}$ must be unitarily equivalent to $%
(S_{\alpha _{0}}\otimes I,I\otimes shift(\beta _{2},\beta _{3},\cdots ))$,
to guarantee the hyponormality of $(T_{1},T_{2})$ it suffices to apply the
Six-point Test (Lemma \ref{joint hypo}) to $\mathbf{k}=(m,n)$, with $m\geq 0$
and $0\leq n\leq 1$.\newline
\textbf{Subcase 1}: $\ \mathbf{k}=(0,0)$. \ Here%
\begin{equation*}
\begin{tabular}{l}
$M(0,0):=\left( 
\begin{array}{cc}
y_{1}^{2}-y_{0}^{2} & \frac{x_{0}^{2}\beta _{0}}{y_{0}}-\beta _{0}y_{0} \\ 
\frac{x_{0}^{2}\beta _{0}}{y_{0}}-\beta _{0}y_{0} & x_{1}^{2}-x_{0}^{2}%
\end{array}%
\right) \geq 0$ \\ 
\\ 
$\Leftrightarrow y_{0}^{2}(y_{1}^{2}-y_{0}^{2})(x_{1}^{2}-x_{0}^{2})\geq
\beta _{0}^{2}(x_{0}^{2}-y_{0}^{2})^{2}$ \\ 
\\ 
$\Leftrightarrow 225\beta _{0}^{2}\leq \frac{35}{72},$%
\end{tabular}%
\end{equation*}%
so $M(0,0)\geq 0$ if and only if 
\begin{equation}
\beta _{0}^{2}\leq \frac{7}{3240}\cong 0.00216.  \label{beta1}
\end{equation}%
\textbf{Subcase 2}: \ $\mathbf{k}=(m,0)$ $($all $m\geq 1)$. \ Fix $m\geq 1$
and let $P_{m}:=\prod_{k=0}^{m-1}x_{k}$ and $Q_{m}:=\prod_{k=0}^{m-1}y_{k}$.
\ We have 
\begin{eqnarray*}
M(m,0) &:&=\left( 
\begin{array}{cc}
y_{m+1}^{2}-y_{m}^{2} & x_{m}\beta _{0}\frac{x_{m}P_{m}}{y_{m}Q_{m}}%
-y_{m}\beta _{0}\frac{P_{m}}{Q_{m}} \\ 
x_{m}\beta _{0}\frac{x_{m}P_{m}}{y_{m}Q_{m}}-y_{m}\beta _{0}\frac{P_{m}}{%
Q_{m}} & \frac{\beta _{1}^{2}\alpha _{0}^{2}}{P_{m}^{2}}-\beta _{0}^{2}\frac{%
P_{m}^{2}}{Q_{m}^{2}}%
\end{array}%
\right) \geq 0 \\
&& \\
&\Leftrightarrow &y_{m}^{2}(y_{m+1}^{2}-y_{m}^{2})(\frac{\beta
_{1}^{2}\alpha _{0}^{2}}{P_{m}^{2}}-\beta _{0}^{2}\frac{P_{m}^{2}}{Q_{m}^{2}}%
)\geq \beta _{0}^{2}(y_{m}^{2}-x_{m}^{2})^{2}\frac{P_{m}^{2}}{Q_{m}^{2}} \\
&& \\
&\Leftrightarrow &\frac{Q_{m}^{2}}{P_{m}^{4}}-\beta _{0}^{2}\geq \frac{%
225(m+2)^{2}(m+3)}{18m+35}\beta _{0}^{2} \\
&& \\
&\Leftrightarrow &\frac{Q_{m}^{2}}{P_{m}^{4}}\geq (\frac{225(m+2)^{2}(m+3)}{%
18m+35}+1)\beta _{0}^{2}.
\end{eqnarray*}%
It follows that $M(m,0)\geq 0$\ \ $($all $m\geq 1)$ if and only if 
\begin{eqnarray*}
\beta _{0}^{2} &\leq &f(m):=\frac{Q_{m}^{2}}{P_{m}^{4}(\frac{%
225(m+2)^{2}(m+3)}{18m+35}+1)} \\
&& \\
&=&\frac{18m+35}{225m^{3}+1575m^{2}+3618m+2735}\prod_{k=0}^{m-1}\frac{%
18k^{2}+71k+70}{9k^{2}+30k+25}\;\;(\text{all }m\geq 1).
\end{eqnarray*}%
Since $f$ is an increasing function of $m$, we see that $M(m,0)\geq 0$\ \ $($%
all $m\geq 1)$ if and only if 
\begin{equation}
\beta _{0}^{2}\leq f(1)=\frac{742}{40765}\cong 0.018.  \label{beta2}
\end{equation}%
\textbf{Subcase 3}: $\ \mathbf{k}=(0,1)$. \ Here we have 
\begin{eqnarray}
M(0,1) &:&=\left( 
\begin{array}{cc}
x_{1}^{2}-x_{0}^{2} & \frac{\alpha _{0}^{2}\beta _{1}}{x_{0}}-\beta _{1}x_{0}
\\ 
\frac{\alpha _{0}^{2}\beta _{1}}{x_{0}}-\beta _{1}x_{0} & \beta
_{2}^{2}-\beta _{1}^{2}%
\end{array}%
\right) \geq 0  \notag \\
&&  \notag \\
&=&\left( 
\begin{array}{cc}
\frac{1}{6} & \frac{\beta _{1}}{x_{0}}(\alpha _{0}^{2}-x_{0}^{2}) \\ 
\frac{\beta _{1}}{x_{0}}(\alpha _{0}^{2}-x_{0}^{2}) & 15\beta _{1}^{2}%
\end{array}%
\right) \geq 0  \notag \\
&&  \notag \\
&\Leftrightarrow &(\alpha _{0}^{2}-\frac{5}{2})^{2}\leq \frac{25}{4},
\label{condition1}
\end{eqnarray}%
which certainly holds, since $0<\alpha _{0}<1$.\newline
\textbf{Subcase 4}: \ $\mathbf{k}=(m,1)$, with $m\geq 1$. \ As in Subcase 2,
fix $m\geq 1$ and let $P_{m}:=\prod_{k=0}^{m-1}x_{k}$. \ We then see that 
\begin{equation}
\begin{tabular}{l}
$M(m,1):=\left( 
\begin{array}{cc}
x_{m+1}^{2}-x_{m}^{2} & \frac{\alpha _{0}\beta _{1}}{x_{m}P_{m}}-x_{m}\frac{%
\alpha _{0}\beta _{1}}{P_{m}} \\ 
\frac{\alpha _{0}\beta _{1}}{x_{m}P_{m}}-x_{m}\frac{\alpha _{0}\beta _{1}}{%
P_{m}} & \beta _{2}^{2}-\frac{\alpha _{0}^{2}\beta _{1}^{2}}{P_{m}^{2}}%
\end{array}%
\right) \geq 0$ \\ 
\\ 
$\Leftrightarrow x_{m}^{2}(x_{m+1}^{2}-x_{m}^{2})(\beta
_{2}^{2}P_{m}^{2}-\alpha _{0}^{2}\beta _{1}^{2})\geq \alpha _{0}^{2}\beta
_{1}^{2}(1-x_{m}^{2})^{2}$ \\ 
\\ 
$\Leftrightarrow x_{m}^{2}(\beta _{2}^{2}P_{m}^{2}-1)\geq
(m+2)(m+3)(1-x_{m}^{2})^{2}$ \\ 
\\ 
$\Leftrightarrow \beta _{2}^{2}\geq g(m):=\frac{1}{P_{m}^{2}}(1+\frac{%
(m+3)(2m+3)^{2}}{3m+5}).$%
\end{tabular}
\label{condition2b}
\end{equation}%
It follows that $M(m,1)\geq 0\;\;$(all $m\geq 1$) if and only if $\beta _{2}$
can be chosen to satisfy (\ref{condition2b}) for all $m\geq 1$. \ Since $g$
is a decreasing function of $m$, it suffices to guarantee that $\beta
_{2}^{2}\geq g(1)=\frac{27}{5}$. \ If we now recall that $\beta _{2}=4\beta
_{1}$ and that $\beta _{1}=\frac{1}{\alpha _{0}}$, this condition is
equivalent to $\alpha _{0}^{2}\leq \frac{80}{27}$, which always holds, since 
$\alpha _{0}<1$. \ 

Therefore, by Subcases 1, 2, 3 and 4, which yield the condition (\ref{beta1}%
), we see that $(T_{1},T_{2})$ is hyponormal if and only if $\beta
_{0}^{2}\leq \frac{7}{3240}$. \ Finally, and since we clearly have $\beta
_{0}<\beta _{1}<\beta _{2}$, we can use the construction in \cite{Sta} to
define $W_{\beta }$, which incidentally has a $2$-atomic Berger measure (cf. %
\cite{RGWSII}).\newline
\textbf{Case 3}: $k_{2}\geq 3$. \ Here we take $p$ and $q$ as in Case 2, to
ensure that the restriction of $\mathbf{T}$ to the subspace associated with
subindices $(m,n)$ with $n\geq k_{2}-2$ is hyponormal. \ Once this is done,
we use Lemma \ref{lemofhypoflat} to obtain $r$, so that the restriction of $%
\mathbf{T}$ to the subspace associated with subindices $(m,n)$ with $n\geq
k_{2}-3$ is hyponormal. \ Repeated application of Lemma \ref{lemofhypoflat}
now completes the proof.
\end{proof}

\begin{corollary}
Theorem \ref{thm33} is optimal.
\end{corollary}

\section{\label{Section4}Propagation in the Subnormal Case}

In this section, we show that Theorem \ref{thm33} can be improved if one of
the $T_{i}$'s is quadratically hyponormal and the other is subnormal. \ In
Theorem \ref{propagation of joint subnormal} and Theorem \ref%
{optimalsubnormal}, we consider horizontal flatness and optimality, and in
Theorem \ref{subnormalflat}, we show that a subnormal $2$-variable weighted
shift with two horizontally consecutive equal weights and two vertically
consecutive equal weights must necessarily be flat. \ As in the previous
section, we then establish that our result is optimal (see Example \ref%
{numericalex1} below). \ We begin with some definitions and preliminary
results.

\begin{definition}
(\cite{CuYo1}) \ Let $\mu $ and $\nu $ be two positive measures on $\mathbb{R%
}_{+}.$ \ We say that $\mu \leq \nu $ on $X:=\mathbb{R}_{+},$ if $\mu
(E)\leq \nu (E)$ for all Borel subset $E\subseteq \mathbb{R}_{+}$;
equivalently, $\mu \leq \nu $ if and only if $\int fd\mu \leq \int fd\nu $
for all $f\in C(X)$ such that $f\geq 0$ on $\mathbb{R}_{+}$.
\end{definition}

\begin{definition}
(\cite{CuYo1}) \ Let $\mu $ be a probability measure on $X\times Y$, and
assume that $\frac{1}{t}\in L^{1}(\mu ).$ \ The extremal measure $\mu _{ext}$
(which is also a probability measure) on $X\times Y$ is given by 
\begin{equation*}
d\mu _{ext}(s,t):=(1-\delta _{0}(t))\frac{1}{t\left\Vert \frac{1}{t}%
\right\Vert _{L^{1}(\mu )}}d\mu (s,t).
\end{equation*}
\end{definition}

\begin{definition}
(\cite{CuYo1})\label{defmarg}Given a measure $\mu $ on $X\times Y$, the
marginal measure $\mu ^{X}$ is given by $\mu ^{X}:=\mu \circ \pi _{X}^{-1}$,
where $\pi _{X}:X\times Y\rightarrow X$ is the canonical projection onto $X$%
. \ Thus, $\mu ^{X}(E)=\mu (E\times Y)$, for every $E\subseteq X$. \ Observe
that if $\mu $ is a probability measure, then so is $\mu ^{X}$.
\end{definition}

\begin{lemma}
(\cite{CuYo1})\label{backward}(Subnormal backward extension of a $1$%
-variable weighted shift) (cf \cite{QHWS}) \ Let $T\equiv shift(\beta
_{0},\beta _{1},\cdots )$ be a unilateral weighted shift whose restriction $%
T_{\mathcal{M}}$ to $\mathcal{M}:=\vee \{e_{1},e_{2},\cdots \}$ is
subnormal, with Berger measure $\eta _{\mathcal{M}}$. \ Then $T$ is
subnormal (with measure $\eta $) if and only if\newline
(i) $\ \frac{1}{t}\in L^{1}(\eta _{\mathcal{M}})$;\newline
(ii) $\ \beta _{0}^{2}\leq (\left\| \frac{1}{t}\right\| _{L^{1}(\eta _{%
\mathcal{M}})})^{-1}$.\newline
In this case, $d\eta (t)=\frac{\beta _{0}^{2}}{t}d\eta _{\mathcal{M}%
}(t)+(1-\beta _{0}^{2}\left\| \frac{1}{t}\right\| _{L^{1}(\eta _{\mathcal{M}%
})})d\delta _{0}(t)$, where $\delta _{0}$ denotes the Dirac measure at $0$.
\ In particular, $T$ is never subnormal when $\eta _{\mathcal{M}}(\{0\})>0$.
\end{lemma}

\begin{lemma}
(\cite{CuYo1})\label{backext}(Subnormal backward extension of a $2$-variable
weighted shift) \ Consider the following $2$-variable weighted shift (see
Figure \ref{Figure 6}), and let $\mathcal{M}$ be the subspace associated to
indices $\mathbf{k}$ with $k_{2}\geq 1$. \ Assume that $\mathbf{T|}_{%
\mathcal{M}}$ is subnormal with measure $\mu _{\mathcal{M}}$ and that $%
W_{0}:=shift(\alpha _{00},\alpha _{10},\cdots )$ is subnormal with measure $%
\xi $. \ Then $\mathbf{T}$ is subnormal if and only if\newline
(i) $\ \frac{1}{t}\in L^{1}(\mu _{\mathcal{M}})$;\newline
(ii) $\ \beta _{00}^{2}\leq (\left\| \frac{1}{t}\right\| _{L^{1}(\mu _{%
\mathcal{M}})})^{-1}$;\newline
(iii) $\ \beta _{00}^{2}\left\| \frac{1}{t}\right\| _{L^{1}(\mu _{\mathcal{M}%
})}(\mu _{\mathcal{M}})_{ext}^{X}\leq \xi $.\newline
Moreover, if $\beta _{00}^{2}\left\| \frac{1}{t}\right\| _{L^{1}(\mu _{%
\mathcal{M}})}=1,$ then $(\mu _{\mathcal{M}})_{ext}^{X}=\xi $. \ In the case
when $\mathbf{T}$ is subnormal, the Berger measure $\mu $ of $\mathbf{T}$ is
given by 
\begin{equation*}
d\mu (s,t)=\beta _{00}^{2}\left\| \frac{1}{t}\right\| _{L^{1}(\mu _{\mathcal{%
M}})}d(\mu _{\mathcal{M}})_{ext}(s,t)+(d\xi (s)-\beta _{00}^{2}\left\| \frac{%
1}{t}\right\| _{L^{1}(\mu _{\mathcal{M}})}d(\mu _{\mathcal{M}%
})_{ext}^{X}(s))d\delta _{0}(t).
\end{equation*}%
\setlength{\unitlength}{1mm} \psset{unit=1mm} 
\begin{figure}[th]
\begin{center}
\begin{picture}(150,138)

\psline{->}(20,20)(130,20)
\psline(20,40)(125,40)
\psline(20,60)(125,60)
\psline(20,80)(125,80)
\psline(20,100)(125,100)
\psline(20,120)(125,120)
\psline{->}(20,20)(20,135)
\psline(40,20)(40,125)
\psline(60,20)(60,125)
\psline(80,20)(80,125)
\psline(100,20)(100,125)
\psline(120,20)(120,125)

\put(11,16){\footnotesize{$(0,0)$}}
\put(35,16){\footnotesize{$(1,0)$}}
\put(55,16){\footnotesize{$(2,0)$}}
\put(78,16){\footnotesize{$\cdots$}}
\put(95,16){\footnotesize{$(m,0)$}}
\put(115,16){\footnotesize{$(m+1,0)$}}

\put(27,21){\footnotesize{$\alpha_{0,0}$}}
\put(47,21){\footnotesize{$\alpha_{1,0}$}}
\put(67,21){\footnotesize{$\alpha_{2,0}$}}
\put(87,21){\footnotesize{$\cdots$}}
\put(107,21){\footnotesize{$\alpha_{m,0}$}}
\put(124,21){\footnotesize{$\alpha_{m+1,0}$}}

\put(27,41){\footnotesize{$\alpha_{0,1}$}}
\put(47,41){\footnotesize{$\alpha_{1,1}$}}
\put(67,41){\footnotesize{$\alpha_{2,1}$}}
\put(87,41){\footnotesize{$\cdots$}}
\put(107,41){\footnotesize{$\alpha_{m,1}$}}
\put(124,41){\footnotesize{$\cdots$}}

\put(27,61){\footnotesize{$\alpha_{0,2}$}}
\put(47,61){\footnotesize{$\alpha_{1,2}$}}
\put(67,61){\footnotesize{$\alpha_{2,2}$}}
\put(87,61){\footnotesize{$\cdots$}}
\put(107,61){\footnotesize{$\alpha_{m,2}$}}
\put(124,61){\footnotesize{$\cdots$}}

\put(27,81){\footnotesize{$\cdots$}}
\put(47,81){\footnotesize{$\cdots$}}
\put(67,81){\footnotesize{$\cdots$}}
\put(87,81){\footnotesize{$\cdots$}}
\put(107,81){\footnotesize{$\cdots$}}
\put(124,81){\footnotesize{$\cdots$}}

\put(27,101){\footnotesize{$\alpha_{0,n}$}}
\put(47,101){\footnotesize{$\alpha_{1,n}$}}
\put(67,101){\footnotesize{$\alpha_{2,n}$}}
\put(87,101){\footnotesize{$\cdots$}}
\put(107,101){\footnotesize{$\alpha_{m,n}$}}
\put(124,101){\footnotesize{$\cdots$}}

\put(27,121){\footnotesize{$\alpha_{0,n+1}$}}
\put(47,121){\footnotesize{$\alpha_{1,n+1}$}}
\put(67,121){\footnotesize{$\alpha_{2,n+1}$}}
\put(87,121){\footnotesize{$\cdots$}}
\put(107,121){\footnotesize{$\alpha_{m,n+1}$}}
\put(124,121){\footnotesize{$\cdots$}}

\psline{->}(70,14)(90,14)
\put(79,10){$\rm{T}_1$}

\put(11,38){\footnotesize{$(0,1)$}}
\put(11,58){\footnotesize{$(0,2)$}}
\put(14,78){\footnotesize{$\vdots$}}
\put(11,98){\footnotesize{$(0,n)$}}
\put(4,118){\footnotesize{$(0,n+1)$}}

\psline{->}(13,70)(13,90)
\put(5,80){$\rm{T}_2$}

\put(20,28){\footnotesize{$\beta_{0,0}$}}
\put(20,48){\footnotesize{$\beta_{0,1}$}}
\put(20,68){\footnotesize{$\beta_{0,2}$}}
\put(22,88){\footnotesize{$\vdots$}}
\put(20,108){\footnotesize{$\beta_{0,n}$}}
\put(20,128){\footnotesize{$\beta_{0,n+1}$}}

\put(40,28){\footnotesize{$\sqrt{\frac{\gamma_{1,1}}{\gamma_{1,0}}}$}}
\put(40,48){\footnotesize{$\sqrt{\frac{\gamma_{1,2}}{\gamma_{1,1}}}$}}
\put(40,68){\footnotesize{$\sqrt{\frac{\gamma_{1,3}}{\gamma_{1,2}}}$}}
\put(42,88){\footnotesize{$\vdots$}}
\put(40,108){\footnotesize{$\sqrt{\frac{\gamma_{1,n+1}}{\gamma_{1,n}}}$}}
\put(42,128){\footnotesize{$\vdots$}}

\put(60,28){\footnotesize{$\sqrt{\frac{\gamma_{2,1}}{\gamma_{2,0}}}$}}
\put(60,48){\footnotesize{$\sqrt{\frac{\gamma_{2,2}}{\gamma_{2,1}}}$}}
\put(60,68){\footnotesize{$\sqrt{\frac{\gamma_{2,3}}{\gamma_{2,2}}}$}}
\put(62,88){\footnotesize{$\vdots$}}
\put(60,108){\footnotesize{$\sqrt{\frac{\gamma_{2,n+1}}{\gamma_{2,n}}}$}}
\put(62,128){\footnotesize{$\vdots$}}

\put(82,28){\footnotesize{$\vdots$}}
\put(82,48){\footnotesize{$\vdots$}}
\put(82,68){\footnotesize{$\vdots$}}
\put(82,88){\footnotesize{$\vdots$}}
\put(82,128){\footnotesize{$\vdots$}}

\put(100,28){\footnotesize{$\sqrt{\frac{\gamma_{m,1}}{\gamma_{m,0}}}$}}
\put(100,48){\footnotesize{$\sqrt{\frac{\gamma_{m,2}}{\gamma_{m,1}}}$}}
\put(100,68){\footnotesize{$\sqrt{\frac{\gamma_{m,3}}{\gamma_{m,2}}}$}}
\put(102,88){\footnotesize{$\vdots$}}
\put(100,108){\footnotesize{$\sqrt{\frac{\gamma_{m,n+1}}{\gamma_{m,n}}}$}}
\put(102,128){\footnotesize{$\vdots$}}

\put(122,28){\footnotesize{$\vdots$}}
\put(122,48){\footnotesize{$\vdots$}}
\put(122,68){\footnotesize{$\vdots$}}
\put(122,88){\footnotesize{$\vdots$}}
\put(122,128){\footnotesize{$\vdots$}}

\end{picture}
\end{center}
\caption{Weight diagram of the 2-variable weighted shift in Lemma \ref%
{backext}}
\label{Figure 6}
\end{figure}
\end{lemma}

\begin{lemma}
\label{backextcor}Let $\mathbf{T}\equiv (T_{1},T_{2})$, let $\mathcal{M}$ be
as in Lemma \ref{backext}, and assume that $\mathbf{T|}_{\mathcal{M}}$ is
subnormal with Berger measure $\mu _{\mathcal{M}}\equiv \delta _{1}\times
\eta .$ \ Assume further that $\beta _{00}^{2}=(\left\| \frac{1}{t}\right\|
_{L^{1}(\mu _{\mathcal{M}})})^{-1}$ and that $W_{\alpha
^{(0)}}:=shift(\alpha _{00},\alpha _{10},\alpha _{20},\cdots )$ is
subnormal. \ Then $\mathbf{T}$ is subnormal if and only if $\alpha _{i0}=1$ $%
($all $i\geq 0)$, that is, $W_{\alpha ^{(0)}}$ must necessarily be the
(unweighted) unilateral shift $U_{+}$.
\end{lemma}

\begin{proof}
Assume first that $\mathbf{T}$ is subnormal. \ Since $d\mu _{\mathcal{M}%
}(s,t)\equiv \delta _{1}(s)d\eta (t)$, we must have%
\begin{eqnarray*}
d(\mu _{\mathcal{M}})_{ext}^{X} &=&((1-\delta _{0}(t))\frac{1}{t\left\Vert 
\frac{1}{t}\right\Vert _{L^{1}(\mu _{\mathcal{M}})}}d\mu _{\mathcal{M}%
}(s,t))^{X} \\
&& \\
&=&d\delta _{1}(s)=d\xi _{\alpha ^{(0)}}(s)\text{ (by Lemma \ref{backext}),}
\end{eqnarray*}%
\ where $\xi _{\alpha ^{(0)}}$ denotes the Berger measure of $W_{\alpha
^{(0)}}$. \ It follows that $W_{\alpha ^{(0)}}=U_{+}$.

Conversely, assume that $W_{\alpha ^{(0)}}=U_{+}$. \ By Lemma \ref{backward}%
, $shift(\beta _{00},\beta _{01},\cdots )$ is subnormal, and we let $\tilde{%
\eta}$ denote its Berger measure. \ If we now let $\mu :=\delta _{1}\times 
\tilde{\eta}$, it easily follows that $\mathbf{T}$ is subnormal with Berger
measure $\mu $.
\end{proof}

\begin{theorem}
\label{propagation of joint subnormal}Let $\mathbf{T}\equiv (T_{1},T_{2})$
be commuting and hyponormal.\newline
(i) \ If $T_{1}$ is quadratically hyponormal, if $T_{2}$ is subnormal, and
if $\alpha _{(k_{1},k_{2})+\mathbf{\varepsilon }_{1}}=\alpha
_{(k_{1},k_{2})} $\ for some $k_{1},k_{2}\geq 0$, then $\mathbf{T}$ is
horizontally flat.\newline
(ii) \ If, instead, $T_{1}$ is subnormal, $T_{2}$ is quadratically
hyponormal, and if $\beta _{(k_{1},k_{2})+\mathbf{\varepsilon }_{2}}=\beta
_{(k_{1},k_{2})}$\ for some $k_{1},k_{2}\geq 0$, then\ $\mathbf{T}$ is
vertically flat.
\end{theorem}

\begin{proof}
Without loss of generality, we only consider the horizontally flat case, and
we further assume $k_{2}=2$, that is, $\alpha _{k_{1},2}=\alpha _{k_{1}+1,2}$
for some $k_{1}\geq 0$. \ By Theorem \ref{thm33} and Proposition \ref%
{quadhypo3}, two equal weights occur at level $3$, i.e., $\alpha
_{k_{1},3}=\alpha _{k_{1}+1,3}$. \ Moreover, for every $\ell \geq 2$ we have 
$\alpha _{k_{1},2}=\alpha _{k_{1},\ell }$ $($all $k_{1}\geq 1)$. \ We now
apply Corollary \ref{corbeta1} to obtain $\beta _{(k_{1},1)}=\beta
_{(k_{1},1)+\mathbf{\varepsilon }_{1}}\;\;$(all $k_{1}\geq 1$). \ By the
commuting property (\ref{commuting}), it follows that%
\begin{equation}
\alpha _{k_{1},2}=\alpha _{k_{1},1}=\alpha _{k_{1}+1,1}\;\;(\text{all }%
k_{1}\geq 1),  \label{alpha1}
\end{equation}%
as desired.
\end{proof}

\begin{corollary}
\label{cor propagation of joint subnormal}Let $\mathbf{T\equiv }%
(T_{1},T_{2}) $ be a subnormal $2$-variable weighted shift.\newline
(i) If $\alpha _{(k_{1},k_{2})+\mathbf{\varepsilon }_{1}}=\alpha
_{(k_{1},k_{2})}$ for some $k_{1},k_{2}\geq 0$, then $\mathbf{T}$ is
horizontally flat.\newline
(ii) If, instead, $\beta _{(k_{1},k_{2})+\mathbf{\varepsilon }_{2}}=\beta
_{(k_{1},k_{2})}$ for some $k_{1},k_{2}\geq 0$, then\ $\mathbf{T}$ is
vertically flat.
\end{corollary}

\begin{proof}
Straightforward from Theorem \ref{propagation of joint subnormal}.
\end{proof}

\begin{remark}
Corollary \ref{cor propagation of joint subnormal} can also be obtained as a
direct consequence of Lemma \ref{backext} and Lemma \ref{backextcor}.
\end{remark}

Theorem \ref{propagation of joint subnormal} is optimal in the following
sense: \ propagation does not necessarily extend either to the left ($0$-th
column) or down ($0$-th level). \ We will actually establish a stronger
result, that is, the optimality of Corollary \ref{cor propagation of joint
subnormal}. \ We first review some basic facts.

\begin{proposition}
(\cite{CuYo2}) \ Let 
\begin{equation}
\alpha _{k}:=\left\{ 
\begin{tabular}{ll}
$\sqrt{\frac{1}{2}},$ & $\text{if }k=0$ \\ 
$\sqrt{\frac{2^{k}+\frac{1}{2}}{2^{k}+1}},$ & $\text{if }k\geq 1$.%
\end{tabular}%
\right.  \label{alpha}
\end{equation}%
Then $W_{\alpha }$ is subnormal with Berger measure $\xi _{\alpha }:=\frac{1%
}{3}\delta _{0}(s)+\frac{1}{3}\delta _{\frac{1}{2}}(s)+\frac{1}{3}\delta
_{1}(s)$.
\end{proposition}

\begin{proposition}
(\cite{CuYo2}) \ Let 
\begin{equation*}
\widehat{\alpha _{k}}:=\left\{ 
\begin{tabular}{ll}
$\sqrt{2},$ & $\text{if }k=0$ \\ 
$\sqrt{\frac{2^{k}+1}{2^{k}+\frac{1}{2}}}\text{,}$ & $\text{if }k\geq 1$%
\end{tabular}%
\right.
\end{equation*}%
then $\prod_{n=0}^{\infty }\widehat{\alpha _{k}}=\sqrt{3}$ . \ (Observe that 
$\widehat{\alpha _{k}}=\frac{1}{\alpha _{k}}$, for $\alpha _{k}$ given by (%
\ref{alpha}).)
\end{proposition}

\begin{theorem}
\label{optimalsubnormal}Consider the weighted shift $\mathbf{T}\equiv
(T_{1},T_{2})$ with weight diagram given by Figure \ref{Figure 9}, where $%
y\leq \frac{1}{\sqrt{3}}$. \ Let $W_{0}:=shift(\alpha _{0},\alpha
_{1},\alpha _{2},\cdots )$, with $\alpha _{k}$ given by (\ref{alpha}). \
Then $\mathbf{T}$ is subnormal.\newline
\setlength{\unitlength}{1mm} \psset{unit=1mm} 
\begin{figure}[th]
\begin{center}
\begin{picture}(140,138)

\psline{->}(20,20)(135,20)
\psline(20,40)(125,40)
\psline(20,60)(125,60)
\psline(20,80)(125,80)
\psline(20,100)(125,100)
\psline(20,120)(125,120)
\psline{->}(20,20)(20,135)
\psline(40,20)(40,125)
\psline(60,20)(60,125)
\psline(80,20)(80,125)
\psline(100,20)(100,125)
\psline(120,20)(120,125)

\put(11,16){\footnotesize{$(0,0)$}}
\put(35,16){\footnotesize{$(1,0)$}}
\put(55,16){\footnotesize{$(2,0)$}}
\put(78,16){\footnotesize{$\cdots$}}
\put(95,16){\footnotesize{$(m,0)$}}
\put(115,16){\footnotesize{$(m+1,0)$}}

\put(27,22){\footnotesize{$\sqrt{\frac{1}{2}}$}}
\put(47,22){\footnotesize{$\sqrt{\frac{5}{6}}$}}
\put(67,22){\footnotesize{$\sqrt{\frac{9}{10}}$}}
\put(87,21){\footnotesize{$\cdots$}}
\put(107,22){\footnotesize{$\sqrt{\frac{2^{m}+\frac{1}{2}}{2^{m}+1}}$}}
\put(124,21){\footnotesize{$\cdots$}}

\put(27,42){\footnotesize{$\sqrt{\frac{1}{2}}$}}
\put(47,41){\footnotesize{$1$}}
\put(67,41){\footnotesize{$1$}}
\put(87,41){\footnotesize{$\cdots$}}
\put(107,41){\footnotesize{$1$}}
\put(124,41){\footnotesize{$\cdots$}}

\put(27,62){\footnotesize{$\sqrt{\frac{1}{2}}$}}
\put(47,61){\footnotesize{$1$}}
\put(67,61){\footnotesize{$1$}}
\put(87,61){\footnotesize{$\cdots$}}
\put(107,61){\footnotesize{$1$}}
\put(124,61){\footnotesize{$\cdots$}}

\put(27,81){\footnotesize{$\cdots$}}
\put(47,81){\footnotesize{$\cdots$}}
\put(67,81){\footnotesize{$\cdots$}}
\put(87,81){\footnotesize{$\cdots$}}
\put(107,81){\footnotesize{$\cdots$}}
\put(124,81){\footnotesize{$\cdots$}}

\put(27,102){\footnotesize{$\sqrt{\frac{1}{2}}$}}
\put(47,101){\footnotesize{$1$}}
\put(67,101){\footnotesize{$1$}}
\put(87,101){\footnotesize{$\cdots$}}
\put(107,101){\footnotesize{$1$}}
\put(124,101){\footnotesize{$\cdots$}}

\put(27,122){\footnotesize{$\sqrt{\frac{1}{2}}$}}
\put(47,121){\footnotesize{$1$}}
\put(67,121){\footnotesize{$1$}}
\put(87,121){\footnotesize{$\cdots$}}
\put(107,121){\footnotesize{$1$}}
\put(124,121){\footnotesize{$\cdots$}}

\psline{->}(70,10)(90,10)
\put(79,6){$\rm{T}_1$}

\put(11,38){\footnotesize{$(0,1)$}}
\put(11,58){\footnotesize{$(0,2)$}}
\put(14,78){\footnotesize{$\vdots$}}
\put(11,98){\footnotesize{$(0,n)$}}
\put(4,118){\footnotesize{$(0,n+1)$}}

\psline{->}(10, 70)(10,90)
\put(5,80){$\rm{T}_2$}

\put(20,28){\footnotesize{$y$}}
\put(20,48){\footnotesize{$\sqrt{\frac{2}{3}}$}}
\put(20,68){\footnotesize{$\sqrt{\frac{3}{4}}$}}
\put(22,88){\footnotesize{$\vdots$}}
\put(20,108){\footnotesize{$\sqrt{\frac{n+1}{n+2}}$}}
\put(22,128){\footnotesize{$\vdots$}}

\put(40,28){\footnotesize{$y$}}
\put(40,48){\footnotesize{$\sqrt{\frac{2}{3}}$}}
\put(40,68){\footnotesize{$\sqrt{\frac{3}{4}}$}}
\put(42,88){\footnotesize{$\vdots$}}
\put(40,108){\footnotesize{$\sqrt{\frac{n+1}{n+2}}$}}
\put(42,128){\footnotesize{$\vdots$}}

\put(60,28){\footnotesize{$\frac{y}{\sqrt{\frac{5}{6}}}$}}
\put(60,48){\footnotesize{$\sqrt{\frac{2}{3}}$}}
\put(60,68){\footnotesize{$\sqrt{\frac{3}{4}}$}}
\put(62,88){\footnotesize{$\vdots$}}
\put(60,108){\footnotesize{$\sqrt{\frac{n+1}{n+2}}$}}
\put(62,128){\footnotesize{$\vdots$}}

\put(82,28){\footnotesize{$\vdots$}}
\put(82,48){\footnotesize{$\vdots$}}
\put(82,68){\footnotesize{$\vdots$}}
\put(82,88){\footnotesize{$\vdots$}}
\put(82,108){\footnotesize{$\vdots$}}
\put(82,128){\footnotesize{$\vdots$}}

\put(100,28){\footnotesize{$\frac{y}{\sqrt{2\gamma_{m-1}}}$}}
\put(100,48){\footnotesize{$\sqrt{\frac{2}{3}}$}}
\put(100,68){\footnotesize{$\sqrt{\frac{3}{4}}$}}
\put(102,88){\footnotesize{$\vdots$}}
\put(100,108){\footnotesize{$\sqrt{\frac{n+1}{n+2}}$}}
\put(102,128){\footnotesize{$\vdots$}}

\put(122,28){\footnotesize{$\vdots$}}
\put(122,48){\footnotesize{$\vdots$}}
\put(122,68){\footnotesize{$\vdots$}}
\put(122,88){\footnotesize{$\vdots$}}
\put(122,108){\footnotesize{$\vdots$}}
\put(122,128){\footnotesize{$\vdots$}}

\end{picture}
\end{center}
\caption{Weight diagram of the $2$-variable weighted shift in Theorem \ref%
{optimalsubnormal}}
\label{Figure 9}
\end{figure}
\end{theorem}

\begin{proof}
To check subnormality, we use Lemma \ref{backext}. \ Since 
\begin{equation*}
\xi _{0}=\frac{1}{3}(\delta _{0}+\delta _{\frac{1}{2}}+\delta _{1})
\end{equation*}%
and 
\begin{equation*}
d\mu _{\mathcal{M}}(s,t)=(d\delta _{0}(s)+d\delta _{1}(s))tdt,
\end{equation*}%
we get 
\begin{equation*}
\beta _{00}^{2}\left\| \frac{1}{t}\right\| _{L^{1}(\mu _{\mathcal{M}})}(\mu
_{\mathcal{M}})_{ext}^{X}=y^{2}(\frac{1}{2}\delta _{0}+\frac{1}{2}\delta
_{1}).\ 
\end{equation*}%
Thus, $y^{2}(\frac{1}{2}\delta _{0}+\frac{1}{2}\delta _{1})\leq \frac{1}{6}%
(\delta _{0}+\delta _{1})\leq \xi _{0}$. \ Lemma \ref{backext} now implies
that $\mathbf{T}$ is subnormal.
\end{proof}

\begin{corollary}
Theorem \ref{propagation of joint subnormal} is optimal.
\end{corollary}

\begin{proof}
Straightforward from Theorem \ref{optimalsubnormal}.
\end{proof}

\begin{theorem}
\label{subnormalflat}Let $\mathbf{T}\equiv (T_{1},T_{2})$ be a subnormal $2$%
-variable weighted shift, and assume that \newline
$\alpha _{(k_{1},k_{2})+\mathbf{\varepsilon }_{1}}=\alpha _{(k_{1},k_{2})}$
and $\beta _{(\ell _{1},\ell _{2})+\mathbf{\varepsilon }_{2}}=\beta _{(\ell
_{1},\ell _{2})}$ for some $k_{1},k_{2},\ell _{1},\ell _{2}\geq 0$. \ Then $%
\mathbf{T}$ is flat.
\end{theorem}

\begin{proof}
Straightforward from Theorem \ref{propagation of joint subnormal}.
\end{proof}

\begin{corollary}
Theorem \ref{subnormalflat} is optimal.
\end{corollary}

\begin{proof}
Straightforward from Example \ref{numericalex1} below.
\end{proof}

\section{\label{symm}Symmetrically flat $2$-variable weighted shifts}

Recall that a $2$-variable weighted shift $\mathbf{T}$ is \textit{flat} if $%
\mathbf{T}$ is horizontally and vertically flat, and \textit{symmetrically
flat} if $\mathbf{T}$ is flat and $\alpha _{11}=\beta _{11}$ (cf. Definition %
\ref{def31}). \ In (\cite[Theorem 2.12]{CuYo1}), we produced an example of a
symmetrically flat, contractive, $2$-variable weighted shift $\mathbf{%
T\equiv }(T_{1},T_{2})$ (that is, $\alpha _{11}=\beta _{11}=1$, and $\left\|
T_{1}\right\| \leq 1$ and $\left\| T_{2}\right\| \leq 1$) with $T_{1},T_{2}$
subnormal, such that $\mathbf{T}$ is hyponormal but not subnormal. \ In this
section, we study the class $\mathcal{SFC}$ of symmetrically flat,
contractive, $2$-variable weighted shifts, with $T_{1}$ and $T_{2}$
subnormal, and we give a complete characterization of hyponormality and
subnormality in $\mathcal{SFC}$; our main result is Corollary \ref%
{hyponotsub}, with a concrete criterion for hyponormality and subnormality.

Symmetrically flat $2$-variable weighted shifts are determined by three main
parts: (i) a subnormal shift in the $0$-th row $(shift(x_{0},x_{1},x_{2}%
\cdots )$, with Berger measure $\xi )$; (ii) a subnormal shift in the $0$-th
column $(shift(y_{0},y_{1},y_{2},\cdots )$, with Berger measure $\eta )$;
and (iii) a positive number $a$ (the $\alpha _{01}$ weight) (cf. Figure \ref%
{flat}). \ By \cite[Theorem 3.3]{CuYo2}, the measures $\xi $ and $\eta $ can
be written as 
\begin{equation}
\begin{tabular}{l}
$\xi \equiv p\delta _{0}+q\delta _{1}+[1-(p+q)]\rho $ \\ 
$\eta \equiv u\delta _{0}+v\delta _{1}+[1-(u+v)]\sigma $,%
\end{tabular}
\label{measureeq}
\end{equation}%
where $0<p,q,u,v<1,$ $p+q\leq 1$, $u+v\leq 1$, and $\rho ,\sigma $ are
probability measures with $\rho (\{0\}\cup \{1\})=\sigma (\{0\}\cup \{1\})=0$%
. \ The following lemma is essential to detect joint hyponormality in the
presence of flatness.

\begin{lemma}
\label{corbackext}(\cite[Theorem 5.2]{CuYo1}) \ Let $\mathbf{T\equiv (}%
T_{1},T_{2})$, let $\mathcal{M}$ be the subspace associated to indices $%
\mathbf{k}$ with $k_{2}\geq 1$, and assume that $\mathbf{T|}_{\mathcal{M}}$
is subnormal with Berger measure $\delta _{1}\times \delta _{1}$. \ Assume
further that $T_{1}$ and $T_{2}$ are contractions, that $W_{0}:=shift(\alpha
_{00},\alpha _{10},\cdots )$ is subnormal with Berger measure $\xi $, and
that $T_{2}$ is subnormal. \ Then $\mathbf{T}$ is subnormal.
\end{lemma}

\begin{remark}
\label{rehyposubnormal} Lemma \ref{corbackext} (together with its proof %
\cite[Theorem 5.2]{CuYo1}) reveals that for the $2$-variable weighted shift
given by Figure \ref{totallyflat}, the hyponormality of $T_{2}$ is
equivalent to the subnormality of $\mathbf{T}$, which in turn is equivalent
to the hyponormality of $\mathbf{T}$.
\end{remark}

\begin{theorem}
\label{flathypo}Let $\mathbf{T\equiv (}T_{1},T_{2})\in \mathcal{SFC}$ be
given by Figure \ref{flat}. \ Then $\mathbf{T}$ is hyponormal if and only if%
\begin{equation}
y_{0}\leq h:=\sqrt{\frac{x_{0}^{2}y_{1}^{2}(x_{1}^{2}-x_{0}^{2})}{%
x_{0}^{2}(x_{1}^{2}-x_{0}^{2})+(a^{2}-x_{0}^{2})^{2}}}.  \label{neweq}
\end{equation}%
\setlength{\unitlength}{1mm} \psset{unit=1mm}
\end{theorem}

\begin{figure}[th]
\begin{center}
\begin{picture}(150,138)

\psline{->}(20,20)(130,20)
\psline(20,40)(125,40)
\psline(20,60)(125,60)
\psline(20,80)(125,80)
\psline(20,100)(125,100)
\psline(20,120)(125,120)
\psline{->}(20,20)(20,135)
\psline(40,20)(40,125)
\psline(60,20)(60,125)
\psline(80,20)(80,125)
\psline(100,20)(100,125)
\psline(120,20)(120,125)

\put(11,16){\footnotesize{$(0,0)$}}
\put(35,16){\footnotesize{$(1,0)$}}
\put(55,16){\footnotesize{$(2,0)$}}
\put(78,16){\footnotesize{$\cdots$}}
\put(95,16){\footnotesize{$(n,0)$}}
\put(115,16){\footnotesize{$(n+1,0)$}}

\put(27,21){\footnotesize{$x_{0}$}}
\put(47,21){\footnotesize{$x_{1}$}}
\put(67,21){\footnotesize{$x_{2}$}}
\put(87,21){\footnotesize{$\cdots$}}
\put(107,21){\footnotesize{$x_{n}$}}
\put(124,21){\footnotesize{$\cdots$}}

\put(27,41){\footnotesize{$a$}}
\put(47,41){\footnotesize{$1$}}
\put(67,41){\footnotesize{$1$}}
\put(87,41){\footnotesize{$\cdots$}}
\put(107,41){\footnotesize{$1$}}
\put(124,41){\footnotesize{$\cdots$}}

\put(24,63){\footnotesize{$\frac{a}{\sqrt{\gamma_{1}(\eta _{1})}}$}}
\put(47,61){\footnotesize{$1$}}
\put(67,61){\footnotesize{$1$}}
\put(87,61){\footnotesize{$\cdots$}}
\put(107,61){\footnotesize{$1$}}
\put(124,61){\footnotesize{$\cdots$}}

\put(27,81){\footnotesize{$\cdots$}}
\put(47,81){\footnotesize{$\cdots$}}
\put(67,81){\footnotesize{$\cdots$}}
\put(87,81){\footnotesize{$\cdots$}}
\put(107,81){\footnotesize{$\cdots$}}
\put(124,81){\footnotesize{$\cdots$}}

\put(23,103){\footnotesize{$\frac{a}{\sqrt{\gamma_{n-1}(\eta _{1})}}$}}
\put(47,101){\footnotesize{$1$}}
\put(67,101){\footnotesize{$1$}}
\put(87,101){\footnotesize{$\cdots$}}
\put(107,101){\footnotesize{$1$}}
\put(124,101){\footnotesize{$\cdots$}}

\put(23,123){\footnotesize{$\frac{a}{\sqrt{\gamma_{n}(\eta _{1})}}$}}
\put(47,121){\footnotesize{$1$}}
\put(67,121){\footnotesize{$1$}}
\put(87,121){\footnotesize{$\cdots$}}
\put(107,121){\footnotesize{$1$}}
\put(124,121){\footnotesize{$\cdots$}}

\psline{->}(70,10)(90,10)
\put(79,6){$\rm{T}_1$}

\put(11,38){\footnotesize{$(0,1)$}}
\put(11,58){\footnotesize{$(0,2)$}}
\put(14,78){\footnotesize{$\vdots$}}
\put(11,98){\footnotesize{$(0,n)$}}
\put(6,118){\footnotesize{$(0,n+1)$}}

\psline{->}(10, 70)(10,90)
\put(5,80){$\rm{T}_2$}

\put(20,28){\footnotesize{$y_{0}$}}
\put(20,48){\footnotesize{$y_{1}$}}
\put(20,68){\footnotesize{$y_{2}$}}
\put(22,88){\footnotesize{$\vdots$}}
\put(20,108){\footnotesize{$y_{n}$}}
\put(22,128){\footnotesize{$\vdots$}}

\put(40,28){\footnotesize{$\frac{ay_{0}}{x_{0}}$}}
\put(40,48){\footnotesize{$1$}}
\put(40,68){\footnotesize{$1$}}
\put(42,88){\footnotesize{$\vdots$}}
\put(40,108){\footnotesize{$1$}}
\put(42,128){\footnotesize{$\vdots$}}

\put(60,28){\footnotesize{$\frac{ay_{0}}{x_{0}x_{1}}$}}
\put(60,48){\footnotesize{$1$}}
\put(60,68){\footnotesize{$1$}}
\put(62,88){\footnotesize{$\vdots$}}
\put(60,108){\footnotesize{$1$}}
\put(62,128){\footnotesize{$\vdots$}}

\put(82,28){\footnotesize{$\vdots$}}
\put(82,48){\footnotesize{$\vdots$}}
\put(82,68){\footnotesize{$\vdots$}}
\put(82,88){\footnotesize{$\vdots$}}
\put(82,108){\footnotesize{$\vdots$}}
\put(82,128){\footnotesize{$\vdots$}}

\put(100,28){\footnotesize{$\frac{ay_{0}}{\sqrt{\gamma_{n}(\xi)}}$}}
\put(100,48){\footnotesize{$1$}}
\put(100,68){\footnotesize{$1$}}
\put(102,88){\footnotesize{$\vdots$}}
\put(100,108){\footnotesize{$1$}}
\put(102,128){\footnotesize{$\vdots$}}

\put(122,28){\footnotesize{$\vdots$}}
\put(122,48){\footnotesize{$\vdots$}}
\put(122,68){\footnotesize{$\vdots$}}
\put(122,88){\footnotesize{$\vdots$}}
\put(122,108){\footnotesize{$\vdots$}}
\put(122,128){\footnotesize{$\vdots$}}

\end{picture}
\end{center}
\caption{Weight diagram of a general symmetrically flat, contractive, $2$%
-variable weighted shift; $\protect\eta _{1}$ denotes the Berger measure of $%
shift(y_{1},y_{2},\cdots )$.}
\label{flat}
\end{figure}

\begin{proof}
By Lemma \ref{corbackext} and Remark \ref{rehyposubnormal}, the subnormality
of $T_{1}$ $($resp. $T_{2})$ implies the subnormality of $\mathbf{T|}_{%
\mathcal{N}}$ $($resp. $\mathbf{T|}_{\mathcal{M}})$, where $\mathcal{N}$
(resp. $\mathcal{M}$) is the subspace associated to indices $\mathbf{k}$
with $k_{1}\geq 1$ (resp. indices $\mathbf{k}$ with $k_{2}\geq 1$). \ Thus,
to verify the hyponormality of $\mathbf{T}$ it suffices to apply the
Six-point Test (Lemma \ref{joint hypo}) to $\mathbf{k}=(0,0)$. \ We have%
\begin{equation}
\begin{tabular}{l}
$\left( 
\begin{array}{cc}
x_{1}^{2}-x_{0}^{2} & \frac{a^{2}y_{0}}{x_{0}}-x_{0}y_{0} \\ 
\frac{a^{2}y_{0}}{x_{0}}-x_{0}y_{0} & y_{1}^{2}-y_{0}^{2}%
\end{array}%
\right) \geq 0$ \\ 
\\ 
$\Leftrightarrow x_{0}^{2}(y_{1}^{2}-y_{0}^{2})(x_{1}^{2}-x_{0}^{2})\geq
y_{0}^{2}(a^{2}-x_{0}^{2})^{2}$ \\ 
\\ 
$\Leftrightarrow y_{0}\leq \sqrt{\frac{%
x_{0}^{2}y_{1}^{2}(x_{1}^{2}-x_{0}^{2})}{%
x_{0}^{2}(x_{1}^{2}-x_{0}^{2})+(a^{2}-x_{0}^{2})^{2}}}=h$.%
\end{tabular}
\label{hyp}
\end{equation}

It follows that $\mathbf{T}$ is hyponormal if and only if $y_{0}\leq h$, as
desired.
\end{proof}

We next consider joint subnormality for $2$-variable weighted shifts in $%
\mathcal{SFC}$. \ We recall Berger's Theorem in the $2$-variable case and
the notion of moment of order $\mathbf{k}$ for a pair $(\alpha ,\beta )$
satisfying (\ref{commuting}). \ Given $\mathbf{k}\in \mathbb{Z}_{+}^{2}$,
the moment of $(\alpha ,\beta )$ of order $\mathbf{k}$ is 
\begin{equation*}
\gamma _{\mathbf{k}}\equiv \gamma _{\mathbf{k}}(\alpha ,\beta ):=\left\{ 
\begin{tabular}{ll}
$1,$ & $\text{if }\mathbf{k}=0$ \\ 
$\alpha _{(0,0)}^{2}\cdot ...\cdot \alpha _{(k_{1}-1,0)}^{2},$ & $\text{if }%
k_{1}\geq 1\text{ and }k_{2}=0$ \\ 
$\beta _{(0,0)}^{2}\cdot ...\cdot \beta _{(0,k_{2}-1)}^{2},$ & $\text{if }%
k_{1}=0\text{ and }k_{2}\geq 1$ \\ 
$\alpha _{(0,0)}^{2}\cdot ...\cdot \alpha _{(k_{1}-1,0)}^{2}\cdot \beta
_{(k_{1},0)}^{2}\cdot ...\cdot \beta _{(k_{1},k_{2}-1)}^{2},$ & $\text{if }%
k_{1}\geq 1\text{ and }k_{2}\geq 1.$%
\end{tabular}%
\right.
\end{equation*}%
We remark that, due to the commutativity condition (\ref{commuting}), $%
\gamma _{\mathbf{k}}$ can be computed using any nondecreasing path from $%
(0,0)$ to $(k_{1},k_{2})$.

\begin{lemma}
\label{Berger}(Berger's Theorem, $2$-variable case) (\cite{JeLu}) $\ $A $2$%
-variable weighted shift $\mathbf{T}\equiv (T_{1},T_{2})$ admits a commuting
normal extension if and only if there is a regular Borel probability measure 
$\mu $ defined on the $2$-dimensional rectangle $R=[0,a_{1}]\times \lbrack
0,a_{2}]$ ($a_{i}:=\left\Vert T_{i}\right\Vert ^{2}$) such that $\gamma _{%
\mathbf{k}}=\iint_{R}\mathbf{t}^{\mathbf{k}}d\mu (\mathbf{t}%
):=\iint_{R}s^{k_{1}}t^{k_{2}}d\mu (s,t)$ \ $($all $\mathbf{k\in }\mathbb{Z}%
_{+}^{2}$).
\end{lemma}

\begin{theorem}
\label{flatsub}Let $\mathbf{T\equiv (}T_{1},T_{2})\in \mathcal{SFC}$ $\ $be
given by Figure \ref{flat}. \ Then $\mathbf{T}$ is subnormal if and only if 
\begin{equation*}
y_{0}\leq s:=\min \left\{ \sqrt{\frac{q}{a^{2}}},\sqrt{\frac{p}{\left\| 
\frac{1}{t}\right\| _{L^{1}(\eta _{1})}-a^{2}}}\right\} \text{.}
\end{equation*}
\end{theorem}

\begin{proof}
Consider the subspaces $\mathcal{M}:=\{\mathbf{k\in }\mathbb{Z}%
_{+}^{2}:k_{2}\geq 1\}$ and $\mathcal{P}:=\{\mathbf{k}\in \mathbb{Z}%
_{+}^{2}:k_{1}\geq 1$ and $k_{2}\geq 1\}$, let $\mathbf{T|}_{\mathcal{M}}$
and $\mathbf{T|}_{\mathcal{P}}$ denote the restrictions of $\mathbf{T}$ to $%
\mathcal{M}$ and $\mathcal{P}$, and let $\eta _{1}$ denote the Berger
measure of $shift(y_{1},y_{2},\cdots )$. \ Since $T_{2}$ is subnormal, and
since $\mathbf{T|}_{\mathcal{P}}$ is the restriction of $\mathbf{T|}_{%
\mathcal{M}}$ to the subspace $\mathcal{P}$, we can apply Lemma \ref{backext}
to $\mathbf{T|}_{\mathcal{M}}$ and the subspace $\mathcal{P}$ (therefore
using as initial data the measures $\delta _{1}\times \delta _{1}$ and $\eta
_{1}$) to show that the subnormality of $T_{2}$ implies $a^{2}\delta
_{1}\leq \eta _{1}$, which in turn gives the subnormality of $\mathbf{T|}_{%
\mathcal{M}}$. \ The Berger measure of $\mathbf{T|}_{\mathcal{M}}$, $\mu _{%
\mathcal{M}}$, is then given by 
\begin{equation}
\mu _{\mathcal{M}}=a^{2}\delta _{1}\times \delta _{1}+\delta _{0}\times
(\eta _{1}-a^{2}\delta _{1}).  \label{mum}
\end{equation}%
Once we know this, we apply Lemma \ref{backext} again, this time to the $2$%
-variable weighted shift $\mathbf{T}$ and the subspace $\mathcal{M}$. \
First, observe that $\left\| \frac{1}{t}\right\| _{L^{1}(\mu _{\mathcal{M}%
})}=\left\| \frac{1}{t}\right\| _{L^{1}(\eta _{1})}$ and from (\ref{mum}) we
have%
\begin{eqnarray*}
d(\mu _{\mathcal{M}})_{ext}(s,t) &=&d(a^{2}\delta _{1}\times \delta
_{1}+\delta _{0}\times (\eta _{1}-a^{2}\delta _{1}))_{ext}(s,t) \\
&& \\
&=&(1-\delta _{0}(t))\frac{1}{t\left\| \frac{1}{t}\right\| _{L^{1}(\mu _{%
\mathcal{M}})}}\{a^{2}d\delta _{1}(s)d\delta _{1}(t) \\
&& \\
&&+d\delta _{0}(s)(d\eta _{1}(t)-a^{2}d\delta _{1}(t))\} \\
&& \\
&=&\frac{1}{\left\| \frac{1}{t}\right\| _{L^{1}(\eta _{1})}}\left\{
a^{2}d\delta _{1}(s)\frac{d\delta _{1}(t)}{t}+d\delta _{0}(s)(\frac{d\eta
_{1}(t)}{t}-a^{2}\frac{d\delta _{1}(t))}{t}\right\}
\end{eqnarray*}%
and therefore%
\begin{eqnarray*}
(\mu _{\mathcal{M}})_{ext}^{X} &=&\frac{1}{\left\| \frac{1}{t}\right\|
_{L^{1}(\eta _{1})}}\left\{ a^{2}\delta _{1}+\delta _{0}(\left\| \frac{1}{t}%
\right\| _{L^{1}(\eta _{1})}-a^{2}\right\} \\
&& \\
&=&\left( 1-\frac{a^{2}}{\left\| \frac{1}{t}\right\| _{L^{1}(\eta _{1})}}%
\right) \delta _{0}+\left( \frac{a^{2}}{\left\| \frac{1}{t}\right\|
_{L^{1}(\eta _{1})}}\right) \delta _{1}.
\end{eqnarray*}%
If we now apply Lemma \ref{backext} and recall (\ref{measureeq}), we see
that the necessary and sufficient condition for $\mathbf{T}$ to be subnormal
is 
\begin{equation*}
y_{0}^{2}\left\| \frac{1}{t}\right\| _{L^{1}(\eta _{1})}\left( (1-\frac{a^{2}%
}{\left\| \frac{1}{t}\right\| _{L^{1}(\eta _{1})}})\delta _{0}+\frac{a^{2}}{%
\left\| \frac{1}{t}\right\| _{L^{1}(\eta _{1})}}\delta _{1}\right) \leq
p\delta _{0}+q\delta _{1}+[1-(p+q)]\rho ,
\end{equation*}%
or equivalently, 
\begin{equation*}
\left\{ 
\begin{array}{ccc}
y_{0}^{2}\left( \left\| \frac{1}{t}\right\| _{L^{1}(\eta _{1})}-a^{2}\right)
& \leq & p \\ 
y_{0}^{2}a^{2} & \leq & q.%
\end{array}%
\right.
\end{equation*}%
It follows at once that $\mathbf{T}$ is subnormal if and only if $y_{0}\leq
s $, as desired.
\end{proof}

We summarize Theorems \ref{flathypo} and \ref{flatsub} as follows.

\begin{corollary}
\label{hyponotsub}The commuting subnormal pair $\mathbf{T\equiv (}%
T_{1},T_{2})$ in Figure \ref{flat} is jointly hyponormal and not subnormal
if and only if 
\begin{equation*}
s<y_{0}\leq h.
\end{equation*}
\end{corollary}

Of course we know that $s\leq h$, but a priori we cannot tell whether the
inequality can be strict. \ We will now exhibit a large collection of $2$%
-variable weighted shifts $\mathbf{T}\in \mathcal{SFC}$ such that $\mathbf{T}
$ is hyponormal but not subnormal; we will do this by describing a
collection of values for $x_{0}$, $x_{1}$, $y_{0}$ and $a$ for which $s<h$.
\ To avoid a trivial case, we shall assume $y_{1}<1$. \ We begin with

\begin{lemma}
\label{lempq}For $\xi \equiv p\delta _{0}+q\delta _{1}+(1-p-q)\rho $ as
above, we have\newline
(i) \ $\int s\;d\xi (s)\geq q$\newline
and\newline
(ii) \ $\int (1-s)\;d\xi (s)\geq p.$\newline
In each case, strict inequality holds if and only if $p+q<1$.

\begin{proof}
Straightforward from the form of $\xi $.
\end{proof}
\end{lemma}

\begin{proposition}
\label{proflat1}Let $x_{0}\sqrt{\frac{1-x_{1}^{2}}{1-x_{0}^{2}}}<a<x_{0}$. \
Then $s<h$.
\end{proposition}

\begin{proof}
We first observe that a straightforward calculation reveals that 
\begin{equation}
P:=x_{0}^{2}x_{1}^{2}+a^{2}-a^{2}x_{0}^{2}-x_{0}^{2}>0  \label{p}
\end{equation}%
whenever $x_{0}\sqrt{\frac{1-x_{1}^{2}}{1-x_{0}^{2}}}<a$. \ Now consider 
\begin{eqnarray}
\frac{h^{2}}{y_{1}^{2}}-\frac{1-x_{0}^{2}}{1-a^{2}} &=&\frac{%
x_{0}^{2}(x_{1}^{2}-x_{0}^{2})}{%
x_{0}^{2}(x_{1}^{2}-x_{0}^{2})+(a^{2}-x_{0}^{2})^{2}}-\frac{1-x_{0}^{2}}{%
1-a^{2}}  \notag \\
&&  \notag \\
&=&\frac{(x_{0}^{2}-a^{2})P}{(1-a^{2})\left[
x_{0}^{2}(x_{1}^{2}-x_{0}^{2})+(a^{2}-x_{0}^{2})^{2}\right] }>0.
\label{hineq}
\end{eqnarray}%
Next, we calculate 
\begin{equation}
1-x_{0}^{2}\equiv \int (1-s)d\xi (s)\geq p\;\;\text{(by Lemma \ref{lempq}%
(ii))}.  \label{pineq}
\end{equation}%
Thirdly, we recall that 
\begin{eqnarray}
1 &=&(\int d\eta _{1}(t))^{2}\leq \int td\eta _{1}(t)\int \frac{1}{t}d\eta
_{1}(t)  \notag \\
&&\text{(using Cauchy-Schwartz in }L^{2}(\eta _{1})\text{)}  \notag \\
&=&y_{1}^{2}\left\| \frac{1}{t}\right\| _{L^{1}(\eta _{1})}<\left\| \frac{1}{%
t}\right\| _{L^{1}(\eta _{1})}.  \label{y1}
\end{eqnarray}%
Finally, we have 
\begin{equation}
\frac{p}{\left( \left\| \frac{1}{t}\right\| _{L^{1}(\eta _{1})}-a^{2}\right) 
}<\frac{p}{\left( 1-a^{2}\right) \left\| \frac{1}{t}\right\| _{L^{1}(\eta
_{1})}}\text{ (since }\left\| \frac{1}{t}\right\| _{L^{1}(\eta _{1})}>1\text{%
).}  \label{vineq}
\end{equation}%
We then have 
\begin{eqnarray*}
s^{2} &\leq &\frac{p}{\left\| \frac{1}{t}\right\| _{L^{1}(\eta _{1})}-a^{2}}<%
\frac{p}{(1-a^{2})\left\| \frac{1}{t}\right\| _{L^{1}(\eta _{1})}}\leq \frac{%
1-x_{0}^{2}}{(1-a^{2})\left\| \frac{1}{t}\right\| _{L^{1}(\eta _{1})}}\;\;%
\text{(by (\ref{pineq}) and (\ref{vineq}))} \\
&& \\
&\leq &\frac{(1-x_{0}^{2})y_{1}^{2}}{1-a^{2}}<h^{2}\;\;\text{(by (\ref{y1})
and (\ref{hineq})),}
\end{eqnarray*}%
as desired.
\end{proof}

\begin{proposition}
\label{proflat2}Let $x_{0}=a$, and assume that $p+q<1$. \ Then $s<h$.
\end{proposition}

\begin{proof}
First, observe that $h=y_{1}$ when $x_{0}=a$; cf. (\ref{neweq}). \ Then 
\begin{eqnarray}
s^{2} &\equiv &\min \left\{ \frac{q}{a^{2}},\frac{p}{\left\| \frac{1}{t}%
\right\| _{L^{1}(\eta _{1})}-a^{2}}\right\}  \notag \\
&&  \notag \\
&<&\frac{1-x_{0}^{2}}{(1-a^{2})\left\| \frac{1}{t}\right\| _{L^{1}(\eta
_{1})}}=y_{1}^{2}\;\;\text{(by Lemma \ref{lempq}, (\ref{pineq}) and (\ref{y1}%
))}  \label{sineq} \\
&=&h^{2}\text{,}  \notag
\end{eqnarray}%
as desired.
\end{proof}

We summarize the above facts in the following result.

\begin{theorem}
\label{flatclass}Let $x_{0}\sqrt{\frac{1-x_{1}^{2}}{1-x_{0}^{2}}}<a\leq
x_{0} $, assume that $p+q<1$, and choose $y_{0}$ in the (nonempty!) interval 
$(s,h] $. \ Then the $2$-variable weighted shift $\mathbf{T}\equiv \mathbf{T}%
(x_{0},x_{1},y_{0},a)$ is hyponormal but not subnormal.
\end{theorem}

We conclude this section by describing a class of numerical examples that
illustrates Theorem \ref{flatclass}. \ Consider the $2$-variable weighted
shift whose weight diagram is given by Figure \ref{numericalex}).

\setlength{\unitlength}{1mm} \psset{unit=1mm} 
\begin{figure}[th]
\begin{center}
\begin{picture}(140,138)

\psline{->}(20,20)(135,20)
\psline(20,40)(125,40)
\psline(20,60)(125,60)
\psline(20,80)(125,80)
\psline(20,100)(125,100)
\psline(20,120)(125,120)
\psline{->}(20,20)(20,135)
\psline(40,20)(40,125)
\psline(60,20)(60,125)
\psline(80,20)(80,125)
\psline(100,20)(100,125)
\psline(120,20)(120,125)

\put(11,16){\footnotesize{$(0,0)$}}
\put(35,16){\footnotesize{$(1,0)$}}
\put(55,16){\footnotesize{$(2,0)$}}
\put(78,16){\footnotesize{$\cdots$}}
\put(95,16){\footnotesize{$(n,0)$}}
\put(115,16){\footnotesize{$(n+1,0)$}}

\put(26,22){\footnotesize{$\sqrt{\frac{1}{2}}$}}
\put(46,22){\footnotesize{$\sqrt{\frac{5}{6}}$}}
\put(66,22){\footnotesize{$\sqrt{\frac{9}{10}}$}}
\put(87,21){\footnotesize{$\cdots$}}
\put(105,22){\footnotesize{$\sqrt{\frac{2^{n}+\frac{1}{2}}{2^{n}+1}}$}}
\put(124,21){\footnotesize{$\cdots$}}

\put(28,41){\footnotesize{$a$}}
\put(47,41){\footnotesize{$1$}}
\put(67,41){\footnotesize{$1$}}
\put(87,41){\footnotesize{$\cdots$}}
\put(107,41){\footnotesize{$1$}}
\put(124,41){\footnotesize{$\cdots$}}

\put(27,62){\footnotesize{$a\frac{3}{2\sqrt{2}}$}}
\put(47,61){\footnotesize{$1$}}
\put(67,61){\footnotesize{$1$}}
\put(87,61){\footnotesize{$\cdots$}}
\put(107,61){\footnotesize{$1$}}
\put(124,61){\footnotesize{$\cdots$}}

\put(27,81){\footnotesize{$\cdots$}}
\put(47,81){\footnotesize{$\cdots$}}
\put(67,81){\footnotesize{$\cdots$}}
\put(87,81){\footnotesize{$\cdots$}}
\put(107,81){\footnotesize{$\cdots$}}
\put(124,81){\footnotesize{$\cdots$}}

\put(24,103){\footnotesize{$\frac{a}{\sqrt{\gamma_{n-1}(\eta_{1})}}$}}
\put(47,101){\footnotesize{$1$}}
\put(67,101){\footnotesize{$1$}}
\put(87,101){\footnotesize{$\cdots$}}
\put(107,101){\footnotesize{$1$}}
\put(124,101){\footnotesize{$\cdots$}}

\put(24,123){\footnotesize{$\frac{a}{\sqrt{\gamma_{n}(\eta_{1})}}$}}
\put(47,121){\footnotesize{$1$}}
\put(67,121){\footnotesize{$1$}}
\put(87,121){\footnotesize{$\cdots$}}
\put(107,121){\footnotesize{$1$}}
\put(124,121){\footnotesize{$\cdots$}}

\psline{->}(70,10)(90,10)
\put(79,6){$\rm{T}_1$}

\put(11,38){\footnotesize{$(0,1)$}}
\put(11,58){\footnotesize{$(0,2)$}}
\put(14,78){\footnotesize{$\vdots$}}
\put(11,98){\footnotesize{$(0,n)$}}
\put(6,118){\footnotesize{$(0,n+1)$}}

\psline{->}(10, 70)(10,90)
\put(5,80){$\rm{T}_2$}

\put(20,28){\footnotesize{$\frac{\sqrt{3}r}{2}$}}
\put(20,48){\footnotesize{$\frac{2\sqrt{2}}{3}$}}
\put(20,68){\footnotesize{$\frac{\sqrt{15}}{4}$}}
\put(22,88){\footnotesize{$\vdots$}}
\put(20,108){\footnotesize{$\frac{\sqrt{(n+1)(n+3)}}{(n+2)}$}}
\put(22,128){\footnotesize{$\vdots$}}

\put(40,30){\footnotesize{$\frac{a\frac{\sqrt{3}r}{2}}{\sqrt{\frac{1}{2}}}$}}
\put(40,48){\footnotesize{$1$}}
\put(40,68){\footnotesize{$1$}}
\put(42,88){\footnotesize{$\vdots$}}
\put(40,108){\footnotesize{$1$}}
\put(42,128){\footnotesize{$\vdots$}}

\put(60,30){\footnotesize{$\frac{a\frac{\sqrt{3}r}{2}}{\sqrt{\frac{1}{2}}\sqrt{\frac{5}{6}}}$}}
\put(60,48){\footnotesize{$1$}}
\put(60,68){\footnotesize{$1$}}
\put(62,88){\footnotesize{$\vdots$}}
\put(60,108){\footnotesize{$1$}}
\put(62,128){\footnotesize{$\vdots$}}

\put(82,28){\footnotesize{$\vdots$}}
\put(82,48){\footnotesize{$\vdots$}}
\put(82,68){\footnotesize{$\vdots$}}
\put(82,88){\footnotesize{$\vdots$}}
\put(82,108){\footnotesize{$\vdots$}}
\put(82,128){\footnotesize{$\vdots$}}

\put(100,30){\footnotesize{$\frac{a\frac{\sqrt{3}r}{2}}{\sqrt{\gamma_{n}(\xi)}}$}}
\put(100,48){\footnotesize{$1$}}
\put(100,68){\footnotesize{$1$}}
\put(102,88){\footnotesize{$\vdots$}}
\put(100,108){\footnotesize{$1$}}
\put(102,128){\footnotesize{$\vdots$}}

\put(122,30){\footnotesize{$\vdots$}}
\put(122,48){\footnotesize{$\vdots$}}
\put(122,68){\footnotesize{$\vdots$}}
\put(122,88){\footnotesize{$\vdots$}}
\put(122,108){\footnotesize{$\vdots$}}
\put(122,128){\footnotesize{$\vdots$}}

\end{picture}
\end{center}
\caption{Weight diagram of the $2$-variable weighted shift in Example \ref%
{numericalex1}}
\label{numericalex}
\end{figure}
To analyze this shift, we will need the following auxiliary results, of
independent interest.

\begin{lemma}
\label{general}(cf. \cite{CLY}) \ For $0<r\leq 1$ let 
\begin{equation}
\beta _{n}(r):=\left\{ 
\begin{tabular}{ll}
$\sqrt{\frac{3}{4}}r,$ & $\text{if }n=0$ \\ 
$\sqrt{\frac{(n+1)(n+3)}{(n+2)^{2}}},$ & $\text{if }n\geq 1$.%
\end{tabular}%
\right.  \label{beta}
\end{equation}%
Then $W_{\beta (r)}$ is subnormal.
\end{lemma}

\begin{proof}
On $\left[ 0,1\right] $, consider the probability measure 
\begin{equation}
d\eta (t):=(1-r^{2})d\delta _{0}(t)+\frac{r^{2}}{2}dt+\frac{r^{2}}{2}d\delta
_{1}(t).\   \label{Bergerofgeneral}
\end{equation}%
For $n\geq 1$ we have 
\begin{eqnarray*}
\gamma _{n}(\beta (r)) &=&r^{2}\frac{3}{2^{2}}\cdot \frac{2\cdot 4}{3^{2}}%
\cdot \frac{3\cdot 5}{4^{2}}\cdot \ldots \cdot \frac{n(n+2)}{(n+1)^{2}} \\
&& \\
&=&\frac{(n+2)r^{2}}{2(n+1)}=\frac{r^{2}}{2}\cdot \frac{1}{n+1}+\frac{r^{2}}{%
2}=\int t^{n}d\eta (t).
\end{eqnarray*}%
Thus, $\eta $ is the Berger measure of $W_{\beta (r)}$, so $W_{\beta (r)}$
is subnormal (all $r\in (0,1]$).
\end{proof}

\begin{lemma}
\label{aux}Let%
\begin{equation*}
\widehat{\beta _{n}}:=\sqrt{\frac{(n+2)^{2}}{(n+3)(n+1)}}\;\;\text{(}n\geq 1%
\text{).}
\end{equation*}%
Then $\prod_{n=1}^{\infty }\widehat{\beta _{n}}=\sqrt{\frac{3}{2}}$. \
(Observe that $\widehat{\beta _{n}}=\frac{1}{\beta _{n}}$ (all $n\geq 1$),
if $\beta _{n}$ is given by (\ref{beta}).)
\end{lemma}

\begin{proof}
Observe that%
\begin{equation*}
\begin{tabular}{l}
$\prod_{n=1}^{k}(\widehat{\beta _{n}})^{2}=\prod_{n=1}^{k}\frac{(n+2)^{2}}{%
(n+3)(n+1)}=\frac{3(k+2)}{2(k+3)}$%
\end{tabular}%
\end{equation*}%
which converge to $\frac{3}{2}$ as $k\rightarrow \infty $.
\end{proof}

\begin{example}
\label{numericalex1}(Illustration of Theorem \ref{flatclass}) \ We first
recall the three assembly parts needed for a $2$-variable weighted shift $T$
to be in $\mathcal{SFC}$: (i) a subnormal shift in the $0$-th row $%
(shift(x_{0},x_{1},x_{2}\cdots )$, with Berger measure $\xi )$; (ii) a
subnormal shift in the $0$-th column $(shift(y_{0},y_{1},y_{2},\cdots )$,
with Berger measure $\eta )$; and (iii) a positive number $a$ (the $\alpha
_{01}$ weight). \ Toward (ii) we shall use the shift in Lemma \ref{general},
with Berger measure given by (\ref{Bergerofgeneral}); toward (i) we shall
use the measure$\text{ }$%
\begin{equation*}
\xi :=\frac{1}{3}(\delta _{0}+\delta _{\frac{1}{2}}+\delta _{1})
\end{equation*}%
on $[0,1]$, so that $p=q=\frac{1}{3}$; finally, toward (iii) we will keep $a$
as a parameter. \ The resulting $2$-variable weighted shift will be denoted $%
\mathbf{T}(a;r)$. \ We will now specify the values of $a$ and $r$ that make $%
\mathbf{T}(a;r)$ contractive, hyponormal, and not subnormal. \ To guarantee
that $\mathbf{T}(a;r)$ is a pair of contractions, and using Lemma \ref{aux},
it is easy to see that we need $a\leq \sqrt{\frac{2}{3}}$. \ Next, we
observe that $x_{0}=\sqrt{\frac{1}{2}}$, $x_{1}=\sqrt{\frac{5}{6}}$, and $%
d\eta _{1}(t)=\frac{2}{3}[tdt+d\delta _{1}(t)]\;\;$($t\in \lbrack 0,1]$), so 
$\left\| \frac{1}{t}\right\| _{L^{1}(\eta _{1})}=\frac{4}{3}$ and $y_{1}=%
\sqrt{\frac{8}{9}}$ . \ Moreover, $x_{0}\sqrt{\frac{1-x_{1}^{2}}{1-x_{0}^{2}}%
}=\sqrt{\frac{1}{6}}$. \ By Theorem \ref{flatclass}, we need to keep $a\in (%
\sqrt{\frac{1}{6}},\sqrt{\frac{1}{2}}]$. \ Thus, for $a\in (\sqrt{\frac{1}{6}%
},\sqrt{\frac{1}{2}}]$ we calculate 
\begin{equation*}
h\equiv \sqrt{\frac{x_{0}^{2}y_{1}^{2}(x_{1}^{2}-x_{0}^{2})}{%
x_{0}^{2}(x_{1}^{2}-x_{0}^{2})+(a^{2}-x_{0}^{2})^{2}}}=\frac{2\sqrt{2}}{3%
\sqrt{1+6(a^{2}-\frac{1}{2})^{2}}}
\end{equation*}%
and 
\begin{equation*}
s\equiv \min \left\{ \sqrt{\frac{q}{a^{2}}},\sqrt{\frac{p}{\left( \left\| 
\frac{1}{t}\right\| _{L^{1}(\eta _{1})}-a^{2}\right) }}\right\} =\min \{%
\frac{1}{\sqrt{3a^{2}}},\sqrt{\frac{1}{4-3a^{2}}}\}=\sqrt{\frac{1}{4-3a^{2}}}%
.
\end{equation*}%
Thus, for $a\in (\sqrt{\frac{1}{6}},\sqrt{\frac{1}{2}}]$ we can then choose $%
y_{0}\equiv \sqrt{\frac{3}{4}}r$ in the interval $(\frac{1}{2\sqrt{1-a^{2}}},%
\frac{2\sqrt{2}}{3\sqrt{1+6(a^{2}-\frac{1}{2})^{2}}}]$ and ensure that $%
\mathbf{T}(a;r)$ is hyponormal and not subnormal (cf. Figure \ref{Figure10}%
). \ \textbf{\qed}

\setlength{\unitlength}{1mm} \psset{unit=1mm}

\begin{figure}[th]
\begin{center}
\begin{picture}(90,75)

\psline{->}(0,10)(0,75)
\put(71,6){$a$}
\put(-10,71){$\sqrt{\frac{3}{4}}r$}
\psline(60,10)(60,11)

\put(-3,5){\footnotesize{$\sqrt{\frac{1}{6}}$}}
\put(57,5){\footnotesize{$\sqrt{\frac{1}{2}}$}}

\put(-7,15){\footnotesize{$0.5$}}
\put(0,15){-}
\put(-7,25){\footnotesize{$0.6$}}
\put(0,25){-}
\put(-7,35){\footnotesize{$0.7$}}
\put(0,35){-}
\put(-7,45){\footnotesize{$0.8$}}
\put(0,45){-}
\put(-7,55){\footnotesize{$0.9$}}
\put(0,55){-}

\put(23,50){\footnotesize{$h$}}
\put(42,27){\footnotesize{$s$}}

\psline{->}(0,10)(75,10)
\pscurve[linewidth=2pt](0,38)(10,42)(20,46)(30,51)(40,55)(50,58)(60,59)

\pscurve[linewidth=1pt](0,18.5)(10,19.5)(20,20.7)(30,22.2)(40,23.9)(50,25.9)(60,28.4)

\put(6,69){\footnotesize{$\mathbf{T}$ is hyponormal but not subnormal 
for $(a,\sqrt{\frac{3}{4}}r)$ in this region}}
\psline[linestyle=dashed,dash=3pt 2pt]{<-}(35,37)(46,66)

\end{picture}
\end{center}
\caption{Graphs of $h$ and $s$ on the interval $[\protect\sqrt{\frac{1}{6}},%
\protect\sqrt{\frac{1}{2}}]$.}
\label{Figure10}
\end{figure}
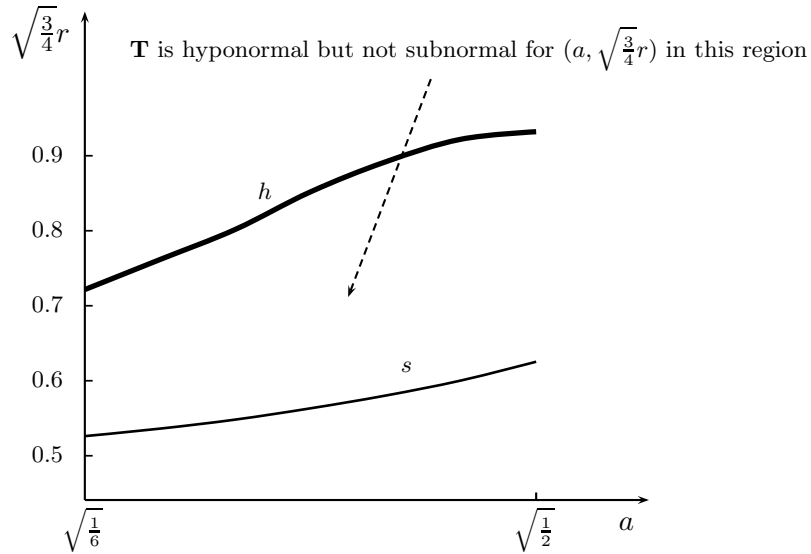
\end{example}

\end{document}